\documentclass[twocolumn,twoside]{IEEEtran}
\usepackage{amsmath,amssymb,amsthm,epsfig,color,subfigure,empheq,graphicx,dsfont}
\usepackage{enumerate,url,algorithm,wasysym,epstopdf,algorithmic,balance}

\DeclareMathOperator{\diag}{dg}

\DeclareMathOperator{\signum}{sign}

\newcommand \tbA{\tilde{\mathbf{A}}}
\newcommand \bA{\mathbf{A}}

\newcommand \bv{\mathbf{v}}

\newcommand \ba{\mathbf{a}}
\newcommand \bp{\mathbf{p}}
\newcommand \bq{\mathbf{q}}
\newcommand \bP{\mathbf{P}}
\newcommand \bQ{\mathbf{Q}}

\newcommand \bX{\mathbf{X}}

\newcommand \bR{\mathbf{R}}
\newcommand \us{\underline{s}}
\newcommand \os{\overline{s}}
\newcommand \ub{\underline{b}}
\newcommand \ob{\overline{b}}

\newtheorem{proposition}{Proposition}
\newtheorem{lemma}{Lemma}
\newtheorem{corollary}{Corollary}
\newtheorem{theorem}{Theorem}

\newtheorem{remark}{Remark}
\allowdisplaybreaks

\begin{document}
\title{Optimal Real-Time Coordination of Energy Storage Units as a Voltage-Constrained Game}

\author{
	Sarthak Gupta, %~\IEEEmembership{Student Member,~IEEE,}
	Vassilis Kekatos,~\IEEEmembership{Senior Member,~IEEE,}
	and
	Walid Saad,~\IEEEmembership{Senior Member,~IEEE}

\thanks{Manuscript submitted October 6, 2017; revised Feb 22, 2018 and Apr 26, 2018; accepted May 27, 2018. Date of publication DATE; date of current version DATE. Paper no. TSG.01451.2017.}

%\thanks{This work was supported by the Spanish Ministry of Education FPU Grant AP2010-1050; CAM Grant S2013/ICE-2933; MINECO Grant TEC2013-41604-R; and NSF grants 1423316, 1442686, 1508993, and 1509040.}

\thanks{The authors are with the Bradley Dept. of ECE, Virginia Tech, Blacksburg, VA 24061, USA. Emails:\{gsarthak,kekatos,walids\}@vt.edu.}

\thanks{Color versions of one or more of the figures is this paper are available online at {http://ieeexplore.ieee.org}.}

\thanks{Digital Object Identifier XXXXXX}
}

\markboth{IEEE TRANSACTIONS ON SMART GRID (to appear)}{Gupta, Kekatos, and Saad: Optimal Real-Time Coordination of Energy Storage Units as a Voltage-Constrained Game}

\maketitle

\begin{abstract}
With increasingly favorable economics and bundling of different grid services, energy storage systems (ESS) are expected to play a key role in integrating renewable generation. This work considers the coordination of ESS owned by customers located at different buses of a distribution grid. Customers participate in frequency regulation and experience energy prices that increase with the total demand. Charging decisions are coupled across time due to battery dynamics, as well as across network nodes due to competitive pricing and voltage regulation constraints. Maximizing the per-user economic benefit while maintaining voltage magnitudes within allowable limits is posed here as a network-constrained game. It is analytically shown that a generalized Nash equilibrium exists and can be expressed as the minimizer of a convex yet infinite-time horizon aggregate optimization problem. To obtain a practical solution, a Lyapunov optimization approach is adopted to design a real-time scheme offering feasible charging decisions with performance guarantees. The proposed method improves over the standard Lyapunov technique via a novel weighting of user costs. By judiciously exploiting the physical grid response, a distributed implementation of the real-time solver is also designed. The features of the novel algorithmic designs are validated using numerical tests on realistic datasets.
\end{abstract}


\section{Introduction}\label{sec:intro}
Energy storage systems (ESS) are expected to lie at the heart of the smart grid, due to their ability to integrate renewables and balance energy \cite{crow-batteries}. Indeed, utility-scale programs for distributed resources and demand response motivate individual customers to employ ESS (including electric vehicle (EV) batteries) for arbitrage, peak shaving, and/or frequency regulation. There is hence a compelling need for ESS control policies to maximize economic benefits while ensuring grid stability. 

Optimal ESS charging schemes can be broadly classified into offline and online. Offline schemes make decisions beforehand by utilizing information about future quantities in the form of exact values and probabilistic or interval characterizations. An offline worst-case ESS coordination scheme is developed in~\cite{YGG13}. Offline protocols for charging EVs under network constraints over a finite horizon are studied in~\cite{ZKG16}. Model predictive control has been also advocated for optimal ESS charging~\cite{Mathieu1}, \cite{Ehsani_mpc}. A stochastic dynamic programming formulation for optimal ESS sizing and control is suggested in~\cite{Harsha}, while \cite{Lohndorf} adopts approximate dynamic programming to jointly store renewable energy and place day-ahead market bids. The underlying assumption about availability of future information renders offline approaches ill-suited for ESS applications with high uncertainty, whereas dynamic programming solutions are impractical for multiple networked ESS. 

Several real-time ESS coordination methods rely on Lyapunov optimization, originally developed for handling data network queues~\cite{Neelybook}. This technique was first adopted for real-time energy arbitrage in data centers in~\cite{Urg11}. Aiming at minimizing the average electricity cost over an infinite time horizon, the derived online scheme yields feasible charging decisions with provable suboptimality guarantees. The technique has been extended to cope with battery charging inefficiencies~\cite{Raja14}; towards distributed ESS implementations involving an aggregator~\cite{SunDongLiang14}; and for handling the exit/return dynamics of EVs~\cite{Sun_EV}. The Lyapunov approach has also been employed for microgrid energy management under network constraints in \cite{Nali17}. However, the feasibility of the charging decisions was only numerically demonstrated. A Lyapunov method for coordinating ESS operating at two timescales is devised in~\cite{G-SIP16}. The Lyapunov technique has been also geared towards managing energy storage over a finite time horizon~\cite{Tli17}. The technique has been interpreted as a stochastic dual approximation algorithm; see \cite{GatsisMarques14} for an application to jointly optimize energy storage and load shedding. Built on competitive analysis, finite time-horizon algorithms have been suggested for arbitrage using energy storage~\cite{Chau16}, and for peak-shaving during electric vehicle charging~\cite{Zhao15}.

The previous approaches presume non-competitive setups. Nevertheless, ESS are usually owned by separate entities and charging decisions are mutually coupled due to competitive pricing or physical network constraints, thus leading naturally to game theoretic formulations. Reference~\cite{Rahi} solves an offline Stackelberg game for ESS charging under behavioral constraints. The competitive scenario in which multiple users aim to minimize their day-ahead cost of operating distributed generation and storage is analyzed in~\cite{Atzeni13}. Sharing ESS resources among users has been shown to be beneficial for arbitrage gains~\cite{Wu16}. Sharing storage and renewable resources have been further studied as coalition games in~\cite{Chis17}.

This work considers the competitive scenario of optimally coordinating user-owned ESS sited at different buses of a distribution grid. The first contribution is to combine Lyapunov optimization with a game-theoretic setup for solving an infinite-horizon energy storage problem. As described in Section~\ref{sec:model}, charging decisions are coupled through voltage constraints and a competitive pricing mechanism incorporating energy charges and frequency regulation benefits. Section~\ref{sec:game} formulates the ESS coordination task as a voltage-constrained game. For this game, we show the existence of a generalized Nash equilibrium that can be found as the minimizer of an aggregate yet infinite-horizon convex quadratic problem. Our second contribution is a \emph{weighted} Lyapunov method to obtain a real-time, near-optimal solver for the aggregate problem (Section~\ref{sec:lyapunov}). Section~\ref{sec:analysis} quantifies the solver's performance gain over the non-weighted formulation of \cite{Tli17}. It is additionally proved that, even under voltage constraints, one can bound the suboptimality and guarantee feasibility of the obtained charging decisions. As a third contribution and to protect customer's privacy, Section~\ref{sec:dist} computes the charging decisions at each control period in a decentralized fashion leveraging dual decomposition and the physical system response. The scheme is numerically tested in Section~\ref{sec:tests}, and Section~\ref{sec:conclusions} concludes the work.

Regarding notation, lower- (upper-) case boldface letters denote column vectors (matrices). Calligraphic symbols are reserved for sets. Vectors $\mathbf{0}$ and $\mathbf{1}$ are the all-zero and all-one vectors. Symbol $\|\mathbf{x}\|_2$ is the $\ell_2$-norm of $\mathbf{x}$ and $^{\top}$ transposition. The main symbols are explained in Table~\ref{tbl:nomenclature}.

%%% NOTATION TABLE
\begin{table}[t]
\renewcommand{\arraystretch}{1.3}
\caption{Nomenclature}
\vspace*{-1em}

\begin{center}
\label{tbl:nomenclature}
\begin{tabular}{|l|l|}
\hline \hline \textbf{Symbol} & \textbf{Meaning}\\
\hline \hline $N$ & number of buses and users\\
\hline $(\ell_n^t,q_n^t)$ & (re)active power demand for user $n$ at time $t$\\
\hline $s_n^t\in [\us_n,\os_n]$ & SoC and its limits for user $n$ at time $t$\\
\hline $b_n^t\in [\ub_n,\ob_n]$ & charge and its limits for user $n$ at time $t$\\
\hline $p_n^t=\ell_n^t+b_n^t$ & net active power demand for user $n$ at time $t$\\
\hline $\tbA$ ($\bA$) & (reduced) branch-bus incidence matrix\\
\hline $I_{m,n}^t$ & current phasor along line $(m,n)$\\
\hline $P_{m,n}^t+jQ_{m,n}^t$ & complex power along line $(m,n)$\\
\hline $r_{m,n}+jx_{m,n}$ & impedance of line $(m,n)$\\
\hline $V_n^t$, $v_n^t$ & voltage phasor and its squared magnitude\\
\hline $\bR,\bX$ & bus resistance and reactance matrices\\
\hline $\alpha,\beta$ & voltage regulation limits\\
\hline $c_0^t\in [\underline{c}_0,\overline{c}_0]$ & base electricity charge and its limits\\
\hline $c_p^t\in [\underline{c}_p,\overline{c}_p]$ & competitive charge and its limits\\
\hline $r^t$ & regulation signal\\
\hline $c_r^t\in [\underline{c}_r,\overline{c}_r]$ & regulation price and its limits\\
\hline $f_n^t(\mathbf{b}^t)$ & electricity cost for user $n$ at time $t$\\
\hline $F_n(\{\mathbf{b}^t\})$ & time-averaged electricity cost for user $n$\\
\hline $f^t(\mathbf{b}^t)$ & aggregate electricity cost at time $t$\\
\hline $F(\{\mathbf{b}^t\})$ & time-averaged aggregate electricity cost\\
\hline $w_n$ & user-specific weight for Lyapunov scheme\\
\hline $\gamma_n$ & SoC-shifting parameter for user $n$\\
\hline $x_n^t=s_n^t+\gamma_n$ & virtual queue for user $n$ at time $t$\\
\hline $\phi,K$ & optimal value and suboptimality gap\\
\hline $x_n^t=s_n^t+\gamma_n$ & virtual queue for user $n$ at time $t$\\
\hline $\delta_n$ & parameter equal to $\frac{\os_n-\us_n+\ub_n-\ob_n}  {\overline{g}_n-\underline{g}_n}$\\
\hline $c_n$ & parameter equal to $w_n x_n- rc_r + c_0$\\
\hline $\underline{g}_n$ & parameter equal to $\underline{c}_0+\frac{\underline{c}_p}{N}\underline{\boldsymbol{\ell}}^\top\mathbf{1} +\frac{\underline{c}_p}{N}\underline{\ell}_n-\overline{c}_r$\\
\hline $\overline{g}_n$ & parameter equal to $\overline{c}_0+\frac{\overline{c}_p}{N}\overline{\boldsymbol{\ell}}^\top\mathbf{1} +\frac{\overline{c}_p}{N}\underline{\ell}_n-\overline{c}_r$\\
\hline $\underline{\lambda}_n,\overline{\lambda}_n$ & Lagrange multipliers for voltage constraints\\
\hline $\eta_\lambda^j$ & step-size at gradient ascent iteration $j$\\
\hline $a$ & total active power demand in dual decomposition\\
\hline \hline
\end{tabular}
\end{center}

\end{table}

%%%%%%%%%%%%%%%%%%%%%%%%%%%%%%%%%%%%%%%%%%%%%%%%%%%%%%%%%%%%%
\section{System Model and Problem Formulation}\label{sec:model}
Consider a power distribution system serving $N$ electricity users indexed by $n$. The system operation is discretized into periods indexed by $t$. Let $\ell_{n}^t$  and $q_{n}^t$ denote respectively the active and reactive load for user $n$ during period $t$. For a compact representation, the loads at period $t$ are collected into the $N$-dimensional vectors $\boldsymbol{\ell}^t$ and $\mathbf{q}^t$. User loads are assumed inelastic and bounded within known intervals as
\begin{subequations}\label{eq:loadlim}
	\begin{align}
\underline{\boldsymbol{\ell}}&\le \boldsymbol{\ell}^t \le\overline{\boldsymbol{\ell}}\label{eq:llim}\\
\underline{\mathbf{q}} &\le \mathbf{q}^t\le\overline{\mathbf{q}}.\label{eq:qlim}
	\end{align}
\end{subequations}

Each user owns an energy storage unit also indexed by $n$. The state of charge (SoC) for unit $n$ at the beginning of slot $t$ is denoted by $s^t_{n}$. The energy by which unit $n$ is charged over period $t$ is denoted by $b^t_{n}$, and it is positive (negative) during (dis)-charging. For simplicity, it is assumed that energy storage units are ideal (unit efficiency). Moreover, since distribution grid customers are currently charged only for active power, energy storage units are assumed to be operated at unit power factor. Upon stacking $\{b^t_{n},s^t_{n}\}_{n=1}^N$ in vectors $(\mathbf{b}^t,\mathbf{s}^t)$ accordingly, the battery dynamics are described as
\begin{subequations}\label{eq:batt}
	\begin{align}
	\mathbf{s}^{t+1}&=\mathbf{s}^t+\mathbf{b}^t \label{eq:soc}\\
	\underline{\mathbf{s}}&\le \mathbf{s}^t\le\overline{\mathbf{s}}\label{eq:slim}\\
	\underline{\mathbf{b}}&\le \mathbf{b}^t\le\overline{\mathbf{b}}\label{eq:blim}
	\end{align}
\end{subequations}
where \eqref{eq:slim} maintains the SoCs within $[\underline{\mathbf{s}},\overline{\mathbf{s}}]$, and \eqref{eq:blim} imposes limits on the charging amount. Customers can reduce their electricity costs by altering their net active power demand $\{p_{n}^t\}_{n=1}^N$ using batteries as
\begin{equation}\label{eq:netpower}
\mathbf{p}^t=\boldsymbol{\ell}^t+\mathbf{b}^t
\end{equation}
where vector $\mathbf{p}^t:=[p_1^t~\cdots~p_N^t]^\top$.

The underlying distribution grid is modeled as a radial single-phase system represented by the graph $\mathcal{G}=(\{0,\mathcal{N}\}, \mathcal{E})$. The substation is indexed by $0$ and the remaining buses comprise the set $\mathcal{N}:=\{1,\dots,N\}$. Each bus hosts one energy storage unit. The edge set $\mathcal{E}$ models distribution lines. If $\pi_n$ is the parent bus of bus $n$, the grid is modeled by the branch flow equations~\cite{BW3}
\begin{subequations}\label{eq:msmv}
	\begin{align}
	-p_n^t&=\sum_{k\in\mathcal{C}_n}P_{n,k}^t  - P_{\pi_n,n}^t + r_{\pi_n,n} |I_{\pi_n,n}^t|^2 \label{eq:mp}\\
	-q_n^t&=\sum_{k\in\mathcal{C}_n}Q_{n,k}^t  - Q_{\pi_n,n}^t + x_{\pi_n,n}|I_{\pi_n,n}^t|^2 \label{eq:mq}\\
	|V_n^t|^2&=|V_{\pi_n}^t|^2 - 2r_{\pi_n,n}P_{\pi_n,n}^t - 2x_{\pi_n,n}Q_{\pi_n,n}^t\nonumber\\
	&\quad+\left(r_{\pi_n,n}^2+x_{\pi_n,n}^2\right)|I_{\pi_n,n}^t|^2 \label{eq:mv}
	\end{align}
\end{subequations}
where $r_{\pi_n,n}+jx_{\pi_n,n}$ is the impedance of line $(\pi_n,n)\in\mathcal{E}$; $\{I_{\pi_n,n}^t,P_{\pi_n,n}^t,Q_{\pi_n,n}^t\}$ is the complex current and the (re)active power flowing from bus $\pi_n$ to bus $n$ at time $t$; $V_n^t$ is the voltage phasor at bus $n$; and $\mathcal{C}_n$ is the set of children buses for $n$. Equations \eqref{eq:mp}--\eqref{eq:mq} stem from conservation of power, while \eqref{eq:mv} captures the drop in squared voltage magnitudes along line $(\pi_n,n)$~\cite{BW3}.

To avoid the nonlinearity in \eqref{eq:mv}, distribution grids are oftentimes studied using the \emph{linear distribution flow} (LDF) model introduced in \cite{BW3}. The latter originates upon setting $I_{\pi_n,n}^t=0$ for all $n$ in \eqref{eq:msmv}. Alternatively, it can be derived using a first-order Taylor series approximation of power injections as functions of voltages evaluated at the flat voltage profile $V_{n}^t=1+j0$ for all $n$~\cite{Deka1}.

The connectivity of a grid with $N$ buses and $L$ distribution lines is captured by the branch-bus incidence matrix $\tbA \in\{0,\pm1\}^{L\times (N+1)}$. Matrix $\tbA$ can be partitioned as $\tbA=[\mathbf{a}_0~\bA]$. For a radial grid with $L=N$, matrix $\bA$ is square and invertible. By dropping the last summands in the RHS of \eqref{eq:msmv}, the LDF model can be compactly expressed as~\cite{VKZG16}
\begin{subequations}\label{eq:mc}
	\begin{align}
	 -\bp^t&=\bA^\top\bP^t\label{eq:mcp}\\
	 -\bq^t&=\bA^\top\bQ^t\label{eq:mcq}\\
	 \bA\bv^t&=2\diag(\{r_{\pi_n,n}\}) \bP^t + 2\diag(\{x_{\pi_n,n}\}) \bQ^t  - \ba_0 v_0\label{eq:mcv}
	\end{align}
\end{subequations}
where $\bv^t:=[|V_1^t|^2~\cdots~|V_N^t|^2]^\top$; the symbol $\diag$ indicates a diagonal matrix; and $v_0=|V_{0}^t|^2$ is the squared voltage magnitude at the substation that is maintained constant. 

Eliminating $(\bP^t,\bQ^t)$ from \eqref{eq:mc} and exploiting the fact that $\mathbf{a}_0 +\bA\mathbf{1}=\mathbf{0}$ or $\bA^{-1}\ba_0=-\mathbf{1}$, the vector $\bv^t$ can be approximated as~\cite{FCL13,Deka1}
\begin{equation}\label{eq:volt}
\bv^t\simeq -\bR\bp^t -\bX\bq^t + v_0\mathbf{1}
\end{equation}
where $\mathbf{1}$ is the all-one vector and
%\begin{subequations}%\label{eq:RX}
\begin{align*}
\bR&:=2\left(\bA^\top \diag^{-1}(\{r_{\pi_n,n}\})\bA\right)^{-1}\\%\label{eq:RX:R}
\bX&:=2\left(\bA^\top \diag^{-1}(\{x_{\pi_n,n}\})\bA\right)^{-1}.%\label{eq:RX:X}
\end{align*}
%\end{subequations}
Numerical tests report that the approximation errors in voltage magnitudes introduced by the LDF model of \eqref{eq:volt} are less than 0.005 pu; see for example~\cite[Fig.~6]{VKZG16}, \cite{BolZamExistence2016}.

Because the entries of $(\mathbf{R},\mathbf{X})$ are non-negative for overhead lines~\cite{VKZG16}, the model in \eqref{eq:volt} implies that voltage magnitudes decrease with increasing $(\mathbf{p}^t,\mathbf{q}^t)$. Grid standards confine nodal voltages to be close to $v_0$~\cite{ansic84}. For example, nodal voltages $|V_n^t|$ should be within $0.97$ and $1.03$ pu, implying that $v_n^t=|V_n^t|^2\in[0.97^2,1.03^2]$ for all $n\in\mathcal{N}$. The latter introduces linear inequality constraints on power demands as
\begin{equation}\label{eq:voltreg} 
\alpha\mathbf{1}\leq-\mathbf{Rp}^t-\mathbf{Xq}^t\leq\beta\mathbf{1}
\end{equation}
where $\alpha:=0.97^2-v_{0,\text{pu}}$ and $\beta:=1.03^2-v_{0,\text{pu}}$. The voltage regulation constraints of \eqref{eq:voltreg} couple charging decisions spatially across costumers. Additional network constraints, such as apparent power limits for lines and (substation) transformers could be included. It is henceforth assumed these limits have been taken care of while allocating loads and sizing transformers, or when the utility grants permission to energy storage installations. For this reason, the next formulations focus on voltage regulation constraints. Nevertheless, if flow limits have to enforced in real time as well, the modification described later in Remark~\ref{re:otherconstraints} can be adopted.

The cost of electricity is varying across time and consists of two components: (i) an energy charge related to real-time energy prices; and (ii) a balancing charge compensating users for participating in frequency regulation. In detail, the cost of electricity for user $n$ at time $t$ is
\begin{equation}\label{eq:cost}
f_n^t(b_n^t,\mathbf{b}_{-n}^t)=\left(c^t_0+ c^t_p\sum_{i=1}^{N}p_{i}^t\right) p_{n}^t-r^tc^t_r b_{n}^t
\end{equation}
where $\mathbf{b}_{-n}^t$ denotes a vector containing the charging decisions for all but the $n$-th user.

The first summand in the right-hand side (RHS) of \eqref{eq:cost} constitutes the energy charge for user $n$. Different from \cite{SunDongLiang14}, the per-unit price is an affine function of the total demand and is assumed to be positive for all $t$: it includes the base charge $c^t_0$ plus the competitive term $c^t_p\sum_{i=1}^{N}p_{i}^t$. When $\sum_{i=1}^{N}p_{i}^t>0$, the per-unit price increases with increasing net total demand and users are motivated to reduce consumption and/or inject energy. When $\sum_{i=1}^{N}p_{i}^t<0$, the per-unit price decreases with increasing energy surplus, thus signaling users to consume. The demand in the feeder may be partially supplied by the transmission grid through the substation.

\begin{remark}\label{re:pricing}
The affine dependence of the electricity price $c^t_0+ c^t_p\sum_{i=1}^{N}p_{i}^t$ on the total demand reflects the fact that the utility participates in a bulk electricity market: higher demand translates to increasingly higher costs. The regulated scenario where customers are subjected to fixed pricing can be captured by setting $c_p^t=0$. Although piecewise-linear pricing could be accommodated, the exposition is restricted to affine pricing to avoid mathematical clutter.
\end{remark}

The second summand in the RHS of \eqref{eq:cost} is the balancing charge defined as the product between the regulation signal $r^t$; the regulation price $c^t_r$; and the battery charge $b^t_n$. The regulation signal is issued by the operator: $r^t=+1$ when there is energy surplus and storage units can only be charged, and $r^t=-1$ during energy deficit periods when storage units can only be discharged. Hence, for all $n$ and $t$
\begin{equation}\label{eq:sign}
r^t=\signum(b^t_n).
\end{equation}
Due to \eqref{eq:sign}, the regulation benefit $r^tc^t_rb^t_n=c_r^t|b^t_n|$ is always positive, and it can thus reduce the total cost for user $n$ in \eqref{eq:cost}.
Here, prices can vary in an arbitrary manner, they are bounded within $0\leq\underline{c}_0\le c^{t}_0\le\overline{c}_0$, $0\leq\underline{c}_p\le c^{t}_p\le\overline{c}_p$ and $0\leq\underline{c}_t\le c^t_{r}\le\overline{c}_r$. The setup where (dis)-charging decisions do not have to comply with $r^t$ as in \eqref{eq:sign} is treated in Remark~\ref{re:no-r}.

%%%%%%%%%%%%%%%%%%%%%%%%%%%%%%%%%%%%%%%%%%%%%%%%%%%%%%%%%%%%%%%%%%
\section{A Game-Theoretic Perspective }\label{sec:game}
Due to the coupling between the users' decisions, minimizing the electricity costs for all users constitutes a voltage-constrained non-cooperative game~\cite{Saad12b}. Each user $n$ seeks to minimize its time-averaged expected electricity cost
\begin{equation}\label{eq:Cn}		
F_n(\{b_{n}^t,\mathbf{b}_{-n}^t\})=\lim_{T\to\infty}\frac{1}{T}\sum_{t=0}^{T-1}\mathbb{E}[f_n^t(b_n^t,\mathbf{b}_{-n}^t)]
\end{equation}
where $\mathbb{E}$ is the expectation over the involved random variables $\{r^t,c^t_0,c^t_{p},c^t_{r},\boldsymbol{\ell}^t,\mathbf{q}^t\}_{t=1}^T$. Then, user $n$ would like to solve the infinite-horizon problem
\begin{subequations}\label{eq:GNEP}
	\begin{align}
	\min_{\{b^t_n,s^t_n\}} ~&F_n(\{b_n^t,\mathbf{b}^t_{-n}\})\\
	\text{s.t. } ~&p^t_n=\ell^t_n+b^t_n\label{eq:singp}\\
      &s^{t+1}_n=s^t_n+b^t_n\label{eq:singsoc}\\
	&\us_n\leq s^t_n \leq \os_n\label{eq:singslim}\\
	&\ub_n\leq b^t_n \leq \ob_n\label{eq:singblim}\\
	&\eqref{eq:voltreg}, \eqref{eq:sign}.
	\end{align}
\end{subequations}
Since the instantaneous costs $f_n^t$ in \eqref{eq:cost} depend on the total demand, the average costs $\{F_n\}_{n=1}^N$ depend on the decisions of all users. The optimal charging decisions are further coupled through the voltage regulation constraints in \eqref{eq:voltreg}, thus rendering \eqref{eq:GNEP} a \emph{generalized Nash game}~\cite{Facchinei07}. 
Formally, we define the game in its strategic form with its set of users $\mathcal{N}$; their costs $\{F_n\}_{n=1}^N$; and the space of feasible (satisfying \eqref{eq:voltreg}) strategies $\mathcal{B}$. The feasible strategies for user $n$ can now be defined as $\mathcal{B}_n(\{\mathbf{b}_{-n}^t\}):=\{\{b_n^t\}: \{b_n^t,\mathbf{b}_{-n}^t\}\in\mathcal{B}\}$.

A sequence of charging decisions $\{\tilde{\mathbf{b}}^t\}$ constitutes a generalized Nash equilibrium (GNE) if it solves \emph{simultaneously} the $N$ coupled minimizations in \eqref{eq:GNEP}. Hence, a GNE is a feasible strategy minimizing the per-user cost as long as the remaining users maintain their strategies, that is for all $n$,
\begin{equation}\label{eq:noregret}
F_n(\{\tilde{b}_n^t,\tilde{\mathbf{b}}^t_{-n}\})\leq F_n(\{\breve{b}_n^t,\tilde{\mathbf{b}}^t_{-n}\}),~\forall \breve{b}_n^t\in \mathcal{B}_n(\{\tilde{\mathbf{b}}_{-n}^t\}).
\end{equation}
A GNE may not necessarily exist. Even if it does, finding it is not always computationally tractable~\cite{Facchinei07}. To prove that a GNE exists for the proposed game and devise algorithms for finding a GNE, we will next transform the set of per-user minimizations in \eqref{eq:GNEP} into a single minimization. The minimizer of this aggregate problem is a GNE for \eqref{eq:GNEP}.

To this end, we first introduce two auxiliary functions. The first function is the aggregate cost at time $t$
\begin{align*}%\label{eq:instcost}
f^t(\mathbf{b}^t)&:=\frac{c^t_p}{2}\left(\sum_{n=1}^{N}p_n^t\right)^2 + \sum_{n=1}^N \left[c^t_0 p_n^t +\frac{c^t_p \left(p^t_n\right)^2}{2}-r^tc^t_r b_{n}^t\right]\nonumber\\
&=c^t_0\mathbf{1}^\top\mathbf{p}^t+\frac{c^t_p}{2} (\mathbf{p}^t)^\top \left(\mathbf{I} + \mathbf{1}\mathbf{1}^\top \right)(\mathbf{p}^t)-r^tc^t_r\mathbf{1}^\top\mathbf{b}^t.
\end{align*}
The function $f^t(\mathbf{b}^t)$ is not the sum of the per-user costs $\{f_n^t(\mathbf{b}^t)\}_{n=1}^N$, but has been constructed so that for all $n$
\begin{subequations}\label{eq:prop1}
\begin{align}
[\nabla f^t(\mathbf{b}^t)]_n&=\frac{\partial f_n^t(\mathbf{b}^t)}{\partial b_n^t}\\
[\nabla^2 f^t(\mathbf{b}^t)]_{n,n}&=\frac{\partial^2 f_n^t(\mathbf{b}^t)}{\partial (b_n^t)^2}.
\end{align}
\end{subequations} 
Using \eqref{eq:prop1} in the second-order Taylor series expansion of the quadratic functions $f^t(\mathbf{b}^t)$ and $f_n^t(\mathbf{b}^t)$ yields the key property
\begin{equation}\label{eq:prop2}
f^t(\mathbf{b}^t)-f^t(\breve{b}_n^t,\mathbf{b}^t_{-n})=f_n^t(\mathbf{b}^t)-f_n^t(\breve{b}_n^t,\mathbf{b}^t_{-n}).
\end{equation}

The second function is the time-averaged aggregate cost
\begin{equation}\label{eq:instcostinf}
F(\{\mathbf{b}^t\}):=\lim_{T\to\infty} \frac{1}{T} \sum_{t=0}^{T-1}\mathbb{E}[f^t(\mathbf{b}^t)].
\end{equation}
Again, the function $F(\{\mathbf{b}^t\})$ is not the sum of $\{F_n(\{\mathbf{b}^t\})\}_{n=1}^N$; but it satisfies
\begin{equation}\label{eq:prop3}
F(\{\mathbf{b}^t\})-F(\{\breve{b}_n^t,\mathbf{b}^t_{-n}\})=F_n(\{\mathbf{b}^t\})-F_n(\{\breve{b}_n^t,\mathbf{b}^t_{-n}\})
\end{equation}
for all $n$. The property in \eqref{eq:prop3} follows easily from \eqref{eq:prop2}. In essence $F(\{\mathbf{b}^t\})$ is the exact potential function for \eqref{eq:GNEP}, which casts the game as a \emph{generalized potential game}~\cite{Facchinei2011}.

Consider next the convex minimization problem
\begin{align}\label{eq:centralized}
\tilde{\phi}:=\min_{\{\mathbf{b}^t\}} ~&~ 	F(\{\mathbf{b}^t\})\\
\text{s.to}~&~\eqref{eq:batt},\eqref{eq:netpower},\eqref{eq:voltreg}, \eqref{eq:sign}.\nonumber
\end{align}
Problem \eqref{eq:centralized} relates to the original problem in \eqref{eq:GNEP} as follows.
		
\begin{proposition}\label{prop:centralized}
The minimizer $\{\tilde{\mathbf{b}}^t\}$ of \eqref{eq:centralized} is a GNE for \eqref{eq:GNEP}.
\end{proposition}

\begin{IEEEproof}Because $f^t(\mathbf{b}^t)$ is quadratic in terms of $\mathbf{b}^t$ with a strictly positive definite Hessian matrix $\frac{c^t_p}{2}\left(\mathbf{I} + \mathbf{1}\mathbf{1}^\top\right)$, it is strictly convex. Strict convexity carries over to $F(\{\mathbf{b}^t\})$. Since the constraints are linear, the optimization in \eqref{eq:centralized} enjoys a unique minimizer $\{\tilde{\mathbf{b}}^t\}$ satisfying
\begin{equation}\label{eq:prop1a}
F(\{\tilde{\mathbf{b}}^t\})<F(\{\breve{b}_n^t,\tilde{\mathbf{b}}^t_{-n}\})
\end{equation} 
for all $\breve{b}_n^t\in \mathcal{B}_n(\{\tilde{\mathbf{b}}_{-n}^t\})$ with $\breve{b}_n^t\neq \tilde{b}_n^t$ and $n \in \mathcal{N}$.

Using \eqref{eq:prop3} in \eqref{eq:prop1a} yields for all $n\in\mathcal{N}$
\begin{equation}\label{eq:prop3b}
%F(\{\tilde{\mathbf{b}}^t\})-F(\{\breve{b}_n^t,\tilde{\mathbf{b}}^t_{-n}\})=
F_n(\{\tilde{\mathbf{b}}^t\})<F_n(\{\breve{b}_n^t,\tilde{\mathbf{b}}^t_{-n}\})
\end{equation}
thus proving that $\{\tilde{\mathbf{b}}^t\}$ is a GNE [cf.~\eqref{eq:noregret}].
\end{IEEEproof}

\begin{remark}
Consider the special case in which each cost $f_n^t$ depends only on $b_n^t$; e.g.,  $c_p^t=0$ in \eqref{eq:cost}. Then, the exact potential function for \eqref{eq:GNEP} can be formulated as the sum of the per-user costs, i.e., $F(\{\mathbf{b}^t\})=\sum_{n=1}^N F_n(\{b_n^t\})$. In this case, a minimizer of \eqref{eq:centralized} is not only a GNE for \eqref{eq:GNEP}, but also its social-welfare solution.
\end{remark}

Proposition~\ref{prop:centralized} asserts that identifying a GNE amounts to solving \eqref{eq:centralized}. Since users lack information on the distribution network, problem \eqref{eq:centralized} can be solved centrally by an aggregator. Albeit convex, the minimization in \eqref{eq:centralized} is challenging: Decisions are coupled over the infinite time horizon via \eqref{eq:soc}--\eqref{eq:slim}, and across grid buses via \eqref{eq:voltreg}. Further, coping with the expected electricity cost requires knowing the joint probability density function of $\{r^t,c^t_0,c^t_p,c^t_r,\boldsymbol{\ell}^t,\mathbf{q}^t\}$. Similar problems are oftentimes tackled through approximate dynamic programming schemes, which are computationally intense~\cite{Powell}. Leveraging Lyapunov-based optimization and dual decomposition, a near-optimal real-time solver is put forth next.

%%%%%%%%%%%%%%%%%%%%%%%%%%%%%%
\section{A Real-Time Solver}\label{sec:lyapunov}
To devise a real-time solver for \eqref{eq:centralized}, consider the problem
\begin{subequations}\label{eq:relaxed}
\begin{align}
\phi':=\min_{\{\mathbf{b}^t\}\in \mathcal{B}} ~&~F(\{\mathbf{b}^t\})\label{eq:relaxed:cost}\\
\text{s.to} ~&~\eqref{eq:blim},\eqref{eq:netpower}, \eqref{eq:voltreg},\eqref{eq:sign}\label{eq:relaxed:con1}\\
~&~ \lim_{T\to\infty} \frac{1}{T}\sum_{t=0}^{T-1}\mathbb{E}\big[b^t_n\big]=0,~\forall n \in \mathcal{N}.\label{eq:relaxed:con2}
\end{align}
\end{subequations}
Problem \eqref{eq:relaxed} is derived from \eqref{eq:centralized} by replacing \eqref{eq:soc}--\eqref{eq:slim} by the constraint \eqref{eq:relaxed:con2} on the expected time-averaged battery charging. In fact, every charging sequence $\{\mathbf{b}^t\}$ complying with \eqref{eq:soc}--\eqref{eq:slim} satisfies also \eqref{eq:relaxed:con2}; see \cite{Urg11}, \cite{Raja14}, or \cite{GatsisMarques14} for a proof. Hence, the minimization in \eqref{eq:relaxed} is a relaxation of the optimization problem in \eqref{eq:centralized}. 

%The SoC dynamics of \eqref{eq:soc} can be rewritten as $s^{T}_n=s^{0}_n+\sum_{t=0}^{T-1}b_n^t$. Applying the expectation operator on the latter, dividing by $T$, and taking the limit of $T$ to infinity implies that for all $n\in\mathcal{N}$
%\begin{equation}
%\lim_{T\to\infty} \frac{1}{T} \mathbb{E}\big[s^{T}_n-s^{0}_n\big]= \lim_{T\to\infty}\frac{1}{T}\sum_{t=0}^{T-1}\mathbb{E}\big[b_n^t\big]. \label{eq:proofrel}
%\end{equation}
%Since $s_n^T$ and $s^{0}_n$ are bounded by \eqref{eq:slim}, the RHS of \eqref{eq:proofrel} goes to zero, thus yielding \eqref{eq:relaxed:con2}. Since \eqref{eq:soc}--\eqref{eq:slim} imply \eqref{eq:relaxed:con2}, problem \eqref{eq:relaxed} is a relaxation of \eqref{eq:centralized}, and hence $\phi'\leq\tilde{\phi}$. 

We next adopt the Lyapunov-based techniques of \cite{Neelybook} to devise a real-time approximate solver for the relaxed problem in \eqref{eq:relaxed}. This solver outputs charging decisions $\{\hat{\mathbf{b}}^t\}$ attaining the objective value $\hat{\phi}:=F(\{\hat{\mathbf{b}}^t\})$. In Section~\ref{sec:analysis}, we will show that $\hat{\phi}$ is $\epsilon$-suboptimal for the relaxed problem in \eqref{eq:relaxed} and that $\{\hat{\mathbf{b}}^t\}$ is feasible not only for \eqref{eq:relaxed}, but also for \eqref{eq:centralized}. This implies that
\begin{equation}\label{eq:phifinal}
\tilde{\phi}\leq\hat{\phi} \leq \phi'+\epsilon \leq\tilde{\phi}+\epsilon.
\end{equation}
In other words, the approximate solver of \eqref{eq:relaxed} achieves bounded suboptimality for \eqref{eq:centralized}. The bounds in \eqref{eq:phifinal} refer to $F(\{\hat{\mathbf{b}}^t\})$ and not to $F_n(\{\hat{\mathbf{b}}^t\})$'s. We will show that the sequence $\{\hat{\mathbf{b}}^t\}$ lies within bounded average distance from $\{\tilde{\mathbf{b}}^t\}$, which is the minimizer  of \eqref{eq:centralized} and the GNEP for \eqref{eq:GNEP}.

To proceed with establishing the previous claims, Lyapunov optimization introduces virtual queues and then stabilizes them to satisfy the average constraint in \eqref{eq:relaxed:con2}~\cite{Neelybook}. For each user $n$, introduce a parameter $\gamma_n$ and define the virtual queue as
\begin{equation}\label{eq:xs}
x^t_n:=s^t_n+\gamma_n.
\end{equation}
Define also the \emph{weighted} Lyapunov function as
\begin{equation}\label{eq:lyapunov}
L^t:=\frac{1}{2}\sum_{n=1}^{N}w_n(x^t_n)^2
\end{equation}
where $\{w_n\}_{n=1}^N$ are positive weights we introduce to handle the heterogeneous capacities and charging rates across energy storage units. Parameters $\{\gamma_n,w_n\}$ are stacked in vectors $\boldsymbol{\gamma}$ and $\mathbf{w}$. Next, we derive upper bounds on the expected differences of successive $L^t$'s given the values of virtual queues collected in vector $\mathbf{x}^t$. These upper bounds will help us later quantify the performance of real-time solvers.

\begin{lemma}\label{le:driftpenalty}
The drift function $\Delta^t:=\mathbb{E}\left[L^{t+1}-L^t|\mathbf{x}^t\right]$ is upper bounded by
\begin{equation}\label{eq:driftpenalty}
\Delta^t\leq \mathbb{E}\left[\sum_{n=1}^{N}w_nx^t_nb^t_n|\mathbf{x}^t\right]+\frac{1}{2}\sum_{n=1}^{N}w_n\max\{\ob_n^2,\ub_n^2\}.
\end{equation}
\end{lemma}

\begin{IEEEproof} Being a shifted version of $s^t_n$, the queue $x^t_n$ evolves similarly to \eqref{eq:soc} as $x^{t+1}_n=x^t_n+b^t_n$. By substituting these queue dynamics in the definition of $\Delta^t$, we get
\begin{align*}
\Delta^t&=\mathbb{E}\left[\frac{1}{2}\sum_{n=1}^{N}w_n\left(2x^t_n b^t_n + {(b^t_n)}^2\right) \vert \mathbf{x}^t\right].
\end{align*}
The bound in \eqref{eq:driftpenalty} follows since \eqref{eq:blim} implies ${(b^t_n)}^2\leq \max\{\ob_n^2,\ub_n^2\}$ for all $n$. 
\end{IEEEproof}

Lyapunov optimization derives an approximate yet real-time solution for \eqref{eq:relaxed} by minimizing the instantaneous cost $f^t(\mathbf{b}^t)$ plus the upper bound on $\Delta^t$ provided by Lemma~\ref{le:driftpenalty}:
\begin{align}\label{eq:realtime}
\hat{\mathbf{b}}^t:=\arg\min_{\mathbf{b}^t}~& ~\sum_{n=1}^N w_n x_n^t b_n^t  + f^t(\mathbf{b}^t)\\
\textrm{s.t.}~&~ \eqref{eq:blim},\eqref{eq:netpower},\eqref{eq:voltreg},\eqref{eq:sign}.\nonumber
\end{align}
Although the average constraint \eqref{eq:relaxed:con2} does not appear in \eqref{eq:realtime}, it is implicitly enforced upon convergence~\cite{Neelybook}. It is worth stressing that \eqref{eq:realtime} depends solely on the current realization of $\mathbf{x}^t$ and $\{r^t,c^t_0,c^t_{p},c^t_{r},\boldsymbol{\ell}^t,\mathbf{q}^t\}$ to find the charging decision $\hat{\mathbf{b}}^t$. It can thus be implemented in real time. The minimization in \eqref{eq:realtime} depends on the parameters $(\boldsymbol{\gamma},\mathbf{w})$, and does not enforce the SoC constraints of \eqref{eq:slim}. By properly designing $(\boldsymbol{\gamma},\mathbf{w})$, the next section optimizes the performance of the real-time solver and guarantees that SoCs remain within limits.

\begin{remark}\label{re:otherconstraints}
Power flows may be restricted by transformer ratings and line thermal limits. These constraints can be readily included in the all the preceding minimizations, that is \eqref{eq:GNEP}, \eqref{eq:centralized}, \eqref{eq:relaxed}, and \eqref{eq:realtime}. In this case, the charging decisions obtained from \eqref{eq:realtime} may not yield realizable SoCs. This issue can be easily resolved if together with the line flow constraints, the SoC constraints $\us_n\leq s^t_n+b^t_n \leq \os_n$ are added to \eqref{eq:realtime}; note that $s^t_n$ is known at time $t$. By doing so, the updated $s_n^{t+1}$'s remain within limits. Although this simple adaptation of \eqref{eq:realtime} yields implementable charging decisions, its performance is not necessarily characterized by the analysis of Section~\ref{sec:analysis}. 
\end{remark}

\begin{remark}\label{re:no-r}
An argument similar to Remark~\ref{re:otherconstraints} holds for constraint \eqref{eq:sign}. If the energy storage units do not have to comply with the (dis)-charging signal $r^t$, problem \eqref{eq:realtime} can be solved upon dropping \eqref{eq:sign} and appending the SoC constraints $\us_n\leq s^t_n+b^t_n \leq \os_n$. In this case, if $r^t=+1$ and $\hat{b}_n^t<0$ at time $t$, user $n$ experiences the regulation penalty of $r^tc^t_r\hat{b}^t_n$, and the suboptimality bound of Section~\ref{sec:analysis} may not hold. 
\end{remark}

%%%%%%%%%%%%%%%%%%%%%%%%%%%%%%%%%%%%%%%%%%%%%%%%%%%%%%%%%%%
\section{Analysis of the Real-Time Solver}\label{sec:analysis}
Next, we show that the charging decisions obtained by the real-time solver of \eqref{eq:realtime} are: (i) feasible for the non-relaxed aggregate problem in \eqref{eq:centralized}; and (ii) within bounded distance both in terms of the optimal cost for \eqref{eq:centralized} and the GNE decisions of \eqref{eq:GNEP}. The analysis extends the results of \cite{Urg11} to the networked ESS setup and relies on two assumptions:

\textit{(a1)} For energy storage unit $n\in\mathcal{N}$, its capacity and charge limits satisfy $\os_n-\us_n>\ob_n-\ub_n$.

\textit{(a2)} In absence of energy storage, the (re)active loads $(\boldsymbol{\ell}^t,\mathbf{q}^t)$ can be served without violating the voltage regulation limits, that is $\alpha\mathbf{1}\leq-\mathbf{R}\boldsymbol{\ell}^t-\mathbf{Xq}^t\leq\beta\mathbf{1}$.

Assumption (a1) essentially excludes fast-charging energy storage units and is commonly adopted in energy storage coordination~\cite{Raja14},~\cite{SunDongLiang14},~\cite{GatsisMarques14}. If $\us_n=0$ and $\ub_n=-\ob_n$, this assumption implies that $2\ob_n<\os_n$, or that it takes more than two periods for an empty battery to be fully charged. This is reasonable if one considers a Tesla supercharger, which can fully charge an EV battery within an hour, participating in a real-time energy market with a control period of 5 or 10 minutes as tested in Section~\ref{sec:tests}. 

Assumption (a2) complies with the assumption that energy storage units do not serve voltage regulation purposes. Excluding energy storage, nodal voltages can be maintained within limits through inverters in solar panels or conventional voltage regulation equipment (regulators, capacitor banks). Although energy storage units do not participate in voltage regulation, they do not incur voltage deviations since the problem in \eqref{eq:realtime} enforces \eqref{eq:voltreg} among its constraints. Albeit useful analytically, Section~\ref{sec:tests} includes tests where (a2) is not met.

The SoCs $\hat{\mathbf{s}}^{t+1}=\hat{\mathbf{s}}^t+\hat{\mathbf{b}}^t$ corresponding to the decisions $\{\hat{\mathbf{b}}^t\}$ obtained from \eqref{eq:realtime} are not explicitly constrained within $[\underline{\mathbf{s}},\overline{\mathbf{s}}]$. This property makes it possible to solve \eqref{eq:realtime} in real time. By properly designing the parameters $(\boldsymbol{\gamma},\mathbf{w})$, the minimizers of \eqref{eq:realtime} will be shown to be feasible for \eqref{eq:centralized}. The next property is the key ingredient to that end and is shown in the appendix.

\begin{theorem}\label{th:charac}
Under (a2), the minimizer $\hat{\mathbf{b}}^t$ of \eqref{eq:realtime} satisfies
\begin{enumerate}[(a)]
\item If $x^t_n+\frac{\underline{g}_n}{w_n}\geq0$ and $r^t>0$, then $\hat b^t_n=0$;
\item If $x^t_n+\frac{\overline{g}_n}{w_n}\leq0$ and $r^t<0$, then $\hat b^t_n=0$;
\end{enumerate}
for all $n\in\mathcal{N}$, and where $\underline{g}_n:=\underline{c}_0+\frac{\underline{c}_p}{N}\underline{\boldsymbol{\ell}}^\top\mathbf{1} +\frac{\underline{c}_p}{N}\underline{\ell}_n-\overline{c}_r$ and $\overline{g}_n:=\overline{c}_0+\frac{\overline{c}_p}{N}\overline{\boldsymbol{\ell}}^\top\mathbf{1} +\frac{\overline{c}_p}{N}\overline{\ell}_n +\overline{c}_r$.
\end{theorem}

Building on Theorem~\ref{th:charac}, the minimizer $\hat{\mathbf{b}}^t$ is next shown to yield feasible states of charge, i.e., $\hat{\mathbf{s}}^t\in [\underline{\mathbf{s}},\overline{\mathbf{s}}]$.

\begin{theorem}\label{th:feasibility}
Under (a1), the minimizer $\hat{\mathbf{b}}^t$ of \eqref{eq:realtime} is also feasible for \eqref{eq:centralized} when $(\boldsymbol{\gamma},\mathbf{w})$ satisfy
\begin{subequations}
\begin{align}
&w_n\delta_n\geq 1 \label{eq:mulim}\\
-\frac{\underline{g}_n}{w_n}+\ob_n-\os_n &\leq \gamma_n \leq -\frac{\overline{g}_n}{w_n}+\ub_n-\us_n \label{eq:gamma}
\end{align}
\end{subequations}
where $\delta_n:= (\os_n-\us_n+\ub_n-\ob_n) / (\overline{g}_n-\underline{g}_n)>0$ for all $n\in\mathcal{N}$.
\end{theorem}

Theorem~\ref{th:feasibility}, which is proved in the appendix, asserts that although the complicating time-coupling constraint $\hat{\mathbf{s}}^t\in [\underline{\mathbf{s}},\overline{\mathbf{s}}]$ has been dropped from \eqref{eq:realtime}, it is actually satisfied by proper parameter tuning. Then, the real-time decisions $\hat{\mathbf{b}}^t$ are feasible for the offline aggregate problem in \eqref{eq:centralized}. Being the minimizer of \eqref{eq:realtime}, the sequence $\hat{\mathbf{b}}^t$ is not necessarily the minimizer of \eqref{eq:centralized}. Nonetheless, it is shown next that $\{\hat{\mathbf{b}}^t\}$ features bounded suboptimality. The ensuing lemma will be needed.

\begin{lemma}[\cite{Neelybook}]\label{th:stat}
If $\{r^t,c^t_0,c^t_{p},c^t_{r},\boldsymbol{\ell}^t,\mathbf{q}^t\}$ are independent and identically distributed (iid) over time, there exists a stationary policy, i.e., a policy selecting ${\grave{\mathbf{b}}^t}$ based only on the current realizations of the involved random variables. This policy further satisfies \eqref{eq:blim}, \eqref{eq:netpower}, \eqref{eq:sign}, \eqref{eq:voltreg}, $\mathbb{E}[\grave{\mathbf{b}}^t]=\mathbf{0}$, and $\mathbb{E}[f^t(\grave{\mathbf{b}}^t)]=\tilde{\phi}$.
\end{lemma}

Using Lemma~\ref{th:stat}, it is shown in the appendix that the average aggregate cost attained by the real-time decisions $\hat{\phi}:=F(\{\hat{\mathbf{b}}^t\})$ satisfies the ensuing suboptimality claim.

\begin{theorem}\label{th:subopt}
If $\{r^t,c^t_0,c^t_{p},c^t_{r},\boldsymbol{\ell}^t,\mathbf{q}^t\}$ are iid, it holds that
\begin{equation}\label{eq:subopt}
\tilde{\phi}\leq\hat{\phi}\leq\tilde{\phi}+K
\end{equation}
where $K:=\frac{1}{2}\sum_{n=1}^{N} w_n \max\{\ob_n^2,\ub_n^2\}$.
\end{theorem}

Due to the quadratic cost, the suboptimality bound in terms of the cost is translated to a suboptimality bound on charging decisions as proved in the appendix.

\begin{theorem}\label{th:l2norm}
Let $\{\hat{\mathbf{b}}^t\}$ be the minimizer of \eqref{eq:realtime}, and $\{\tilde{\mathbf{b}}^t\}$ the minimizer of \eqref{eq:centralized} that is also the sought GNE. Then,
\begin{equation*}
\lim_{T\to\infty} \frac{1}{T}\sum_{t=0}^{T-1}\mathbb{E}[\|\hat{\mathbf{b}}^t-{\tilde{\mathbf{b}}}^t\|^2_2]\leq \frac{2K}{\underline{c}_p}.
\end{equation*}
\end{theorem}

Theorem~\ref{th:l2norm} guarantees that the obtained charging decisions lie close to the GNE decisions, thus providing a sense of satisfaction among users. Based on the suboptimality bounds provided by Theorems~\ref{th:subopt} and \ref{th:l2norm}, the performance of the real-time solver can be optimized by minimizing the quantity $K$ over the weights $\mathbf{w}$ subject to $w_n\delta_n\geq 1$ for all $n\in\mathcal{N}$~[cf.~\eqref{eq:mulim}]. Since $K$ is separable over $\{w_n\}$, the optimal weights are simply $w_n^\star:=\delta_n^{-1}$. Moreover, by plugging $\{w_n^\star\}$ into \eqref{eq:gamma}, it is not hard to verify that its leftmost and rightmost sides coincide. Then, the allowable range for each $\gamma_n$ collapses to a single value.

\begin{corollary}\label{co:maxopt}
To minimize the suboptimality bound and guarantee feasibility of the SoC variables, the parameters $\{\boldsymbol{\gamma},\mathbf{w}\}$ in \eqref{eq:realtime} should be selected for all $n\in\mathcal{N}$ as
\begin{subequations}\label{eq:optparams}
\begin{align}
w_n^\star&:=\delta_n^{-1}\label{eq:optparams:w}\\
\gamma_n^\star&:=-\frac{\overline{g}_n(\os_n-\ob_n)-\underline{g}_n(\us_n-\ub_n)}{\overline{g}_n-\underline{g}_n}>0.\label{eq:optparams:gamma}
\end{align}
\end{subequations}
The suboptimality bound becomes $K^\star:=\sum_{n=1}^{N} \frac{\max\{\ob_n^2,\ub_n^2\}}{2\delta_n}$.
\end{corollary}

Corollary~\ref{co:maxopt} sets the values for parameters $\{\boldsymbol{\gamma},\mathbf{w}\}$ in \eqref{eq:realtime}. Since $\os_n-\ob_n>\us_n-\ub_n$ by assumption (a1), it follows that $\gamma_n^\star>\us_n-\ub_n\geq 0$. Corollary~\ref{co:maxopt} further justifies having user-specific weights $\{w_n\}$ in \eqref{eq:realtime}: The standard non-weighted Lyapunov technique of \cite{Urg11} would have resulted in a common weight for all $n$
\begin{equation}\label{eq:wnw}
w_\text{nw}=\delta_{\min}^{-1}
\end{equation}
where $\delta_{\min}:=\min \delta_n$.\footnote{In fact, instead of weighting the first summand in the cost of \eqref{eq:realtime} by $w_\text{nw}$, \cite{Urg11} would equivalently multiply its second summand by $\mu=w_\text{nw}^{-1}=\delta_{\min}$.} The weight $w_\text{nw}$ guarantees there exist $\{\gamma_n\}$ satisfying \eqref{eq:gamma} and attaining suboptimality gap
\begin{equation}\label{eq:Kprime}
K':=\frac{1}{2\delta_{\min}}\sum_{n=1}^{N} \max\{\ob_n^2,\ub_n^2\}\geq K^\star.
\end{equation}
Rather than having the user with the smallest $\delta_n$ controlling the algorithm performance, the weights in \eqref{eq:realtime} account for heterogeneity across energy storage units.

%%%%%%%%%%%%%%%%%%%%%%%%%%%%%%%%%%%%%%%%%%%%%%%%%%%%%%%%%%%%%%%%%%%%%%%%%%%%%%%%%%%%%
\section{Distributed Implementation}\label{sec:dist}
The minimization in \eqref{eq:realtime} can be performed in a centralized fashion using standard (e.g., interior point-based) or customized solvers for linearly-constrained convex quadratic programs. In that case however, the limits $(\underline{\mathbf{s}},\overline{\mathbf{s}},\underline{\mathbf{b}},\overline{\mathbf{b}})$ along with the sequences $\{\boldsymbol{\ell}^t,\mathbf{q}^t,\mathbf{s}^t\}$ need to be communicated from the users to the aggregator. To waive possible concerns on user privacy, a distributed scheme for tackling \eqref{eq:realtime} is proposed next. To simplify notation, the superscript $t$ will be dropped. 

Let us first rewrite \eqref{eq:realtime} in the equivalent form 
\begin{subequations}\label{eq:realtime2}
\begin{align}
\min_{\mathbf{b},a}~&~\sum_{n=1}^N \left[c_n b_n +\frac{c_p}{2}(b_n+\ell_n)^2\right]+\frac{c_p}{2}a^2\label{eq:realtime2:cost}\\
\textrm{s.to}~&~\frac{(1-r)\ub_n}{2}\leq b_n\leq \frac{(1+r)\ob_n}{2},~\forall n\label{eq:realtime2:con1}\\
~&~a=\mathbf{1}^\top (\mathbf{b}+\boldsymbol{\ell})\label{eq:realtime2:con2}\\
~&~\alpha\mathbf{1}\leq-\mathbf{R}(\mathbf{b}+\boldsymbol{\ell})-\mathbf{Xq}\leq\beta\mathbf{1}\label{eq:realtime2:con3}
\end{align}
\end{subequations}
where $c_n:=w_n x_n- rc_r + c_0$. Note that variable $\mathbf{p}=\mathbf{b}+\boldsymbol{\ell}$ has been eliminated; the constraint \eqref{eq:realtime2:con1} combines \eqref{eq:blim} and \eqref{eq:sign}; the new variable $a$ captures the net active power demand through \eqref{eq:realtime2:con2}; and \eqref{eq:realtime2:con3} enforces voltage regulation. 

To derive a decentralized solver, we adopt dual decomposition and introduce Lagrange multipliers $\nu$ and  $\underline{\boldsymbol{\lambda}}\geq \mathbf{0}$ ($\overline{\boldsymbol{\lambda}}\geq \mathbf{0}$) for constraint \eqref{eq:realtime2:con2} and the left-hand (right-hand) side of \eqref{eq:realtime2:con3}, respectively. Dual decomposition updates these Lagrange multipliers through the projected gradient ascent iterations~\cite[Ch.~6]{Be99}
\begin{subequations}\label{eq:dual}
\begin{align}
\nu^{j+1}&=\nu^j+\eta_\nu^j \left[a^j - \mathbf{1}^\top(\mathbf{b}^j +  \boldsymbol{\ell} )\right]\label{eq:dual:c1}\\
\underline{\boldsymbol{\lambda}}^{j+1}&=\max\left\{\underline{\boldsymbol{\lambda}}^j +  \eta_\lambda^j \left[\mathbf{R} (\mathbf{b}^j + \boldsymbol{\ell}) + \mathbf{Xq} + \alpha\mathbf{1} \right], \mathbf{0}\right\}\label{eq:dual:c2}\\
\overline{\boldsymbol{\lambda}}^{j+1}&= \max\left\{\overline{\boldsymbol{\lambda}}^j - \eta_\lambda^j \left[\mathbf{R}(\mathbf{b}^j + \boldsymbol{\ell}) + \mathbf{Xq} + \beta\mathbf{1}\right], \mathbf{0}\right\}\label{eq:dual:c3}
\end{align}
\end{subequations}
where $\eta_\nu^j,\eta_\lambda^j>0$ are step sizes; the maximum operator is applied entrywise; and $(a^j,\mathbf{b}^j)$ are the minimizers of the Lagrangian function associated with the minimization in \eqref{eq:realtime2} evaluated at $(\nu^j,\underline{\boldsymbol{\lambda}}^{j},\overline{\boldsymbol{\lambda}}^j)$. The primal variable $a^j$ can be found by the aggregator in closed-form as
\begin{equation}\label{eq:aprimal}
a^j=\arg\min_{a}\left\{\frac{c_p}{2} a^2+\nu^j a \right\}=-\frac{\nu^j}{c_p}.
\end{equation}
Let $\tilde{\boldsymbol{\lambda}}^j:=\mathbf{R}(\underline{\boldsymbol{\lambda}}^j-\overline{\boldsymbol{\lambda}}^j) - \nu^j\mathbf{1}$. Then, the charging decisions at iteration $j$ can be updated \emph{separately} over users by solving 
\begin{align}\label{eq:user1}
b_n^j:=\arg\min_{b_n}~&~ \frac{c_p}{2} (b_n +\ell_n)^2 +(c_n+\tilde{\lambda}_n^j)b_n\\
\text{s.to}~&~\left(\frac{1-r}{2}\right)\ub_n\leq b_n\leq \left(\frac{1+r}{2}\right)\ob_n.\nonumber
\end{align}
The minimizer of \eqref{eq:user1} can be readily found in closed form as
\begin{equation}\label{eq:user2}
b_n^j=\left[-\frac{c_n+\tilde{\lambda}_n^j}{c_p}-\ell_n\right]_{(1-r)\ub_n/2}^{(1+r)\ob_n/2}
\end{equation}
where the $[x]_{a}^b:=\max\{\min\{x,b\},a\}$ projects $x$ onto the interval $[a,b]$. Given the strict convexity of the objective in \eqref{eq:realtime}, the iterations in \eqref{eq:dual}--\eqref{eq:user2} are guaranteed to converge to the optimal dual and primal variables~\cite{Be99}. The steps involved for solving \eqref{eq:realtime} at time $t$ are tabulated as Algorithm~\ref{alg:1}.

\begin{algorithm}[t]
	\caption{Distributed solver for \eqref{eq:realtime} at time $t$.} \label{alg:1} 
	\begin{algorithmic}[1]
		\STATE Aggregator initializes dual variables $(\nu^{t,0},\underline{\boldsymbol{\lambda}}^{t,0},\overline{\boldsymbol{\lambda}}^{t,0})$.
		\STATE Aggregator observes $\{c_0^t,c_p^t,c_r^t\}$.
		\STATE Aggregator estimates $\mathbf{R}\boldsymbol{\ell}^t+\mathbf{Xq}^t$ as $\mathbf{v}^t-v_0\mathbf{1}-\mathbf{R}\hat{\mathbf{b}}^{t-1}$.
		\STATE Each user $n$ observes $\{r^t,c_0^t,c_p^t,c_r^t,\ell_n^t\}$.
		\FOR{$j=0,1,\ldots,$}
		\STATE Aggregator communicates to users the entries of\\ $\tilde{\boldsymbol{\lambda}}^{t,j}:=\mathbf{R}(\underline{\boldsymbol{\lambda}}^{t,j}-\overline{\boldsymbol{\lambda}}^{t,j}) - \nu^{t,j}\mathbf{1}$.
		\STATE User $n$ updates $b_n^{t,j}$ via \eqref{eq:user2} and communicates it back to the aggregator. 
		\STATE Net load $\mathbf{1}^\top\boldsymbol{\ell}$ is communicated to the aggregator.
		\STATE Aggregator updates primal variable $a^{t,j}$ from \eqref{eq:aprimal}.
		\STATE Aggregator updates multipliers $(\nu^{t,j},\underline{\boldsymbol{\lambda}}^{t,j},\overline{\boldsymbol{\lambda}}^{t,j})$ by \eqref{eq:dual}.
		\ENDFOR
	\end{algorithmic}
\end{algorithm}

To update the dual variables in \eqref{eq:dual}, the aggregator needs to know $\mathbf{R}\boldsymbol{\ell}^t+\mathbf{Xq}^t$ and $\mathbf{1}^\top\boldsymbol{\ell}^t$. Given \eqref{eq:volt}, the former can be calculated indirectly as $\mathbf{R}\boldsymbol{\ell}^t+\mathbf{Xq}^t = \mathbf{v}^t-v_0\mathbf{1}-\mathbf{R}\hat{\mathbf{b}}^{t-1}$ assuming the previous charging decisions $\hat{\mathbf{b}}^{t-1}$ persist at the beginning of period $t$ and that the aggregator measures $\mathbf{v}^t$. Calculating $\mathbf{1}^\top\boldsymbol{\ell}^t$ can also be performed without sending private information to the aggregator. Instead, the total load $\mathbf{1}^\top\boldsymbol{\ell}^t$ can be computed by having nodes communicating over a spanning tree rooted at the aggregator. Leaf nodes pass their load values to their parents, their parents sum up the received information and their own load, and the recursion proceeds. The communication tree does not necessarily match the electric grid and can be randomized at each time. 

Due to the way optimal weights $w_n$'s are determined in Remark~\ref{co:maxopt}, the users do not need to communicate their $\delta_n$ or $w_n$ to the aggregator. This is an added advantage of our \emph{weighted} Lyapunov optimization method over the conventional one where $\delta_n$'s have to be shared to identify $\delta_{\min}$.

%%%%%%%%%%%%%%%%%%%%%%%%%%%%%%%%%%%%%%%%%%%%%%%%%%%%%%%%%%%%%%%%%%%%%%%%%%%%%%%%%%%%%%%%%%
\section{Numerical Tests}\label{sec:tests}
To recapitulate, each user would ideally like to reach the generalized Nash equilibrium (GNE) obtained by solving \eqref{eq:GNEP}. Given its stochastic and infinite-time horizon nature, solving \eqref{eq:GNEP}, or even obtaining its optimal value $\tilde{\phi}$ is computationally intractable. Theorem~\ref{th:l2norm} though ensured that the decisions obtained by \eqref{eq:realtime} lie close to the GNE with the distance being proportional to the suboptimality gap $(\hat{\phi}-\tilde{\phi})$. This section evaluates different real-time charging schemes based on the objective value $\phi:=F(\{\mathbf{b}^t\})$ they attain. This is because if a sequence $\{\mathbf{b}^t\}$ yields smaller $\phi$, this sequence lies closer to the sought GNE $\{\tilde{\mathbf{b}}^t\}$, which is impossible to compute.

The developed charging schemes were evaluated using load data from the Pecan Street project comprising both consumption and solar generation~\cite{pecandata}. Lacking reactive
injections, a lagging power factor of 0.9 was assumed. Five minute averages were obtained from the minute-based load data. The so obtained load data were placed on the IEEE 13-bus and 34-bus feeders along with energy storage units. Both feeders were converted to single-phase grids as described in~\cite{tac2015gan}. Voltage deviations were allowed to lie within $\pm1\%$ by setting $\alpha=-0.0199$ and $\beta=0.020$ in \eqref{eq:voltreg}. 
	 	
The developed Lyapunov-based algorithm was compared against two competing alternatives. The first alternative is a standard Lyapunov-based algorithm. Both the developed and the standard (non-weighted) Lyapunov schemes were operated for the parameter values minimizing the related suboptimality gaps [cf.~\eqref{eq:optparams} and \eqref{eq:wnw}]. The second alternative is the \emph{greedy} charging scheme
\begin{align}\label{eq:greedy}
\min_{\mathbf{b}^t}~& ~f^t(\mathbf{b}^t)\\
\textrm{s.t.}~&~ \eqref{eq:batt},\eqref{eq:netpower},\eqref{eq:voltreg},\eqref{eq:sign}\nonumber
\end{align}
which can be implemented in real time similar to \eqref{eq:realtime}. Different from \eqref{eq:realtime} though, the problem in \eqref{eq:greedy} involves only the instantaneous cost $f^t$, and it explicitly enforces the SoC constraints \eqref{eq:soc}--\eqref{eq:slim}. Since by using $\{\boldsymbol{\gamma}^\star,\mathbf{w}^\star\}$, the minimizer of \eqref{eq:realtime} also satisfies \eqref{eq:soc}--\eqref{eq:slim}, the only difference between \eqref{eq:realtime} and \eqref{eq:greedy} ends up being their costs. When $s_n^t$ is large, the extra term $w_n x_n^t b_n^t$ in the cost of \eqref{eq:realtime} becomes larger (recall $x_n^t=s_n^t+\gamma_n$) hence promoting smaller charging amounts $\hat{b}_n^t$ for the same values of $\{c_0^t,c_r^t,c_p^t,r^t\}$. In other words, the greedy scheme selects the currently optimal decision, whereas \eqref{eq:realtime} takes into account the current price along with the current SoC. For example, the scheme of \eqref{eq:realtime} requires higher financial benefit to decide to charge an almost full battery. 

\begin{figure}[t]
\centering
\includegraphics[scale=0.09]{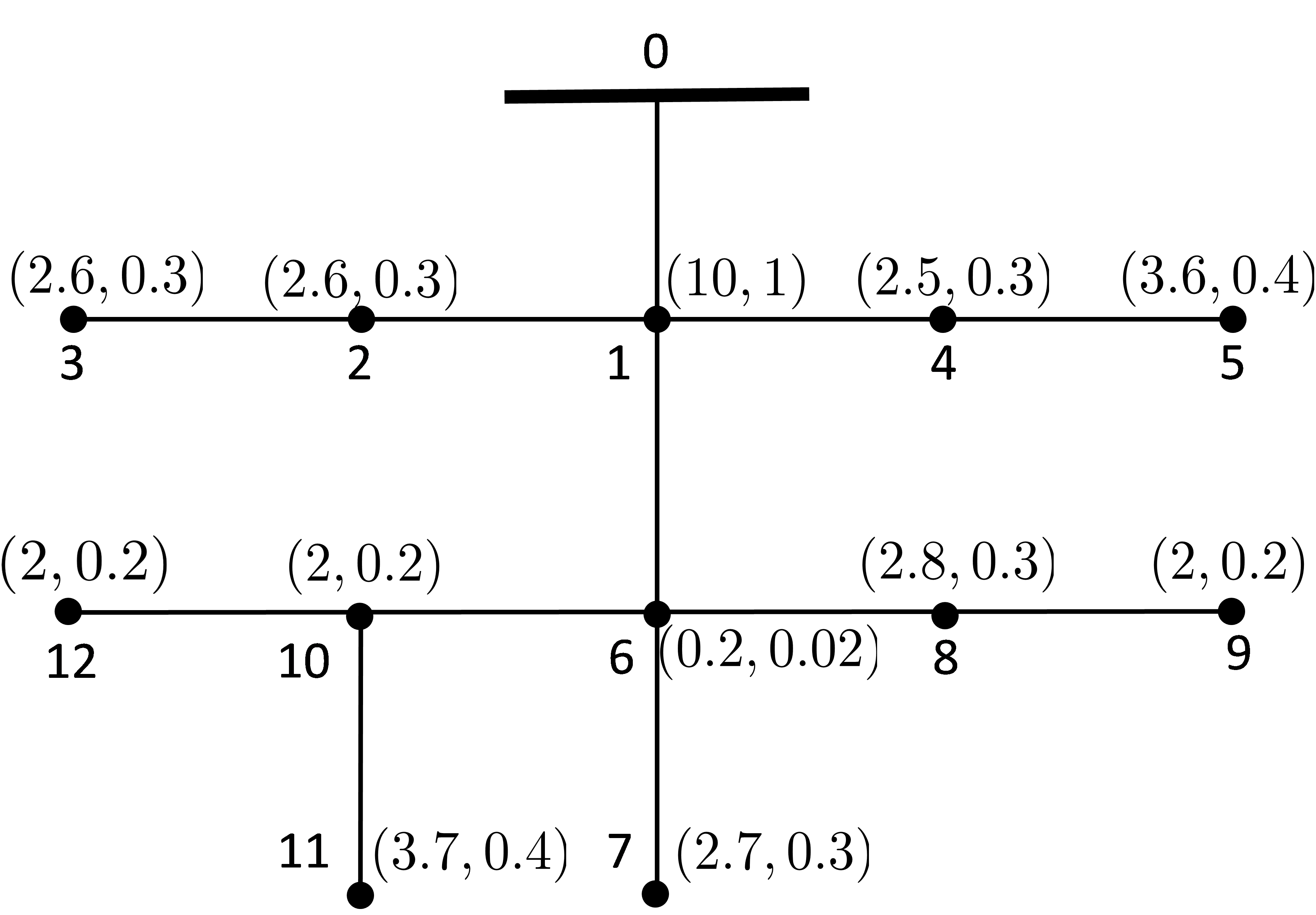}
\caption{IEEE 13-bus feeder for Scenario 1: battery capacities and charging rates are indicated as $(\os_n,\ob_n)$ in $10^{-1}$ and $10^{-2}$ kWh, accordingly.}
%\vspace*{-1em}
\label{fig:ieee13}
\end{figure}

\begin{figure}[t]
\centering
\includegraphics[scale=0.4]{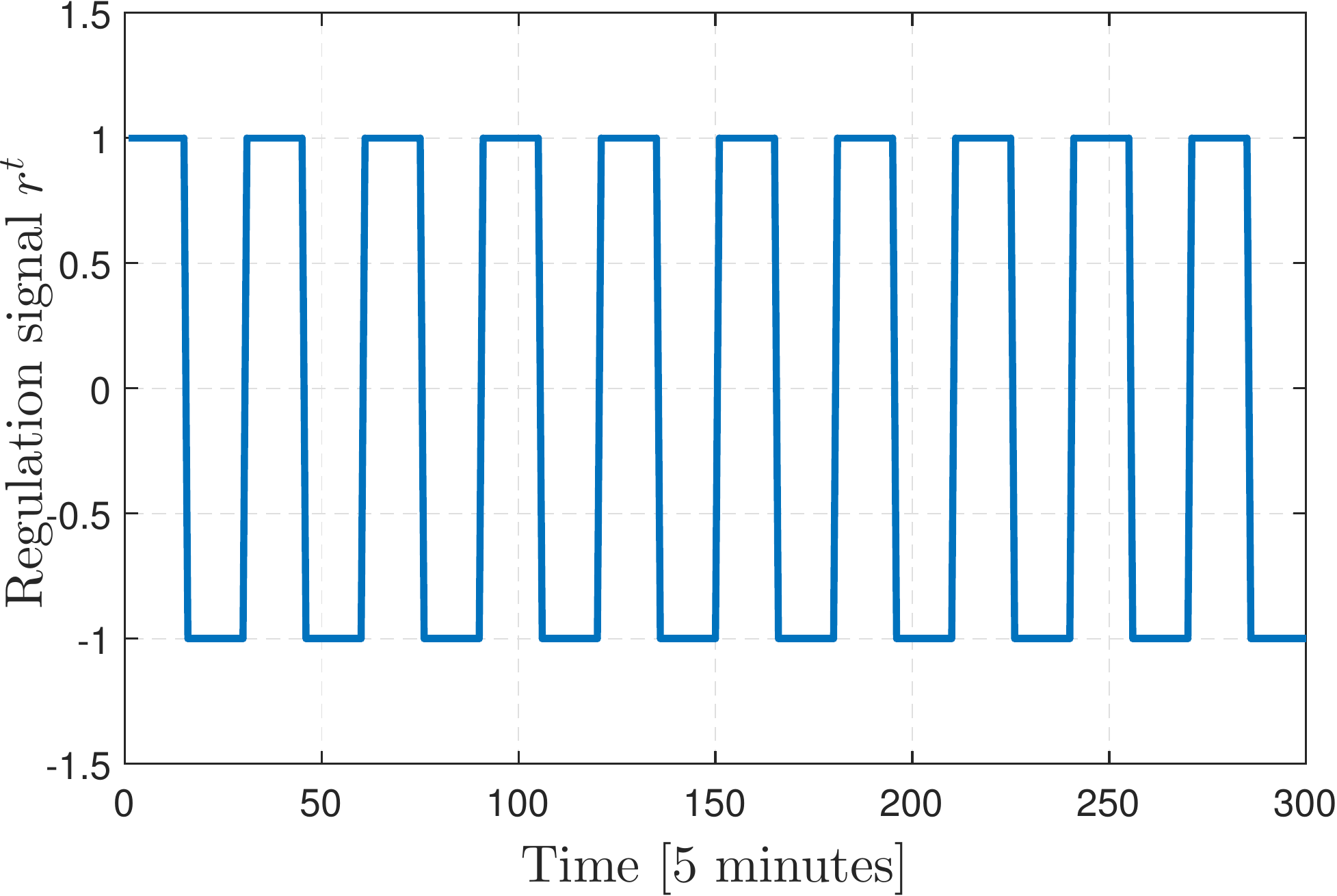}\\
\includegraphics[scale=0.4]{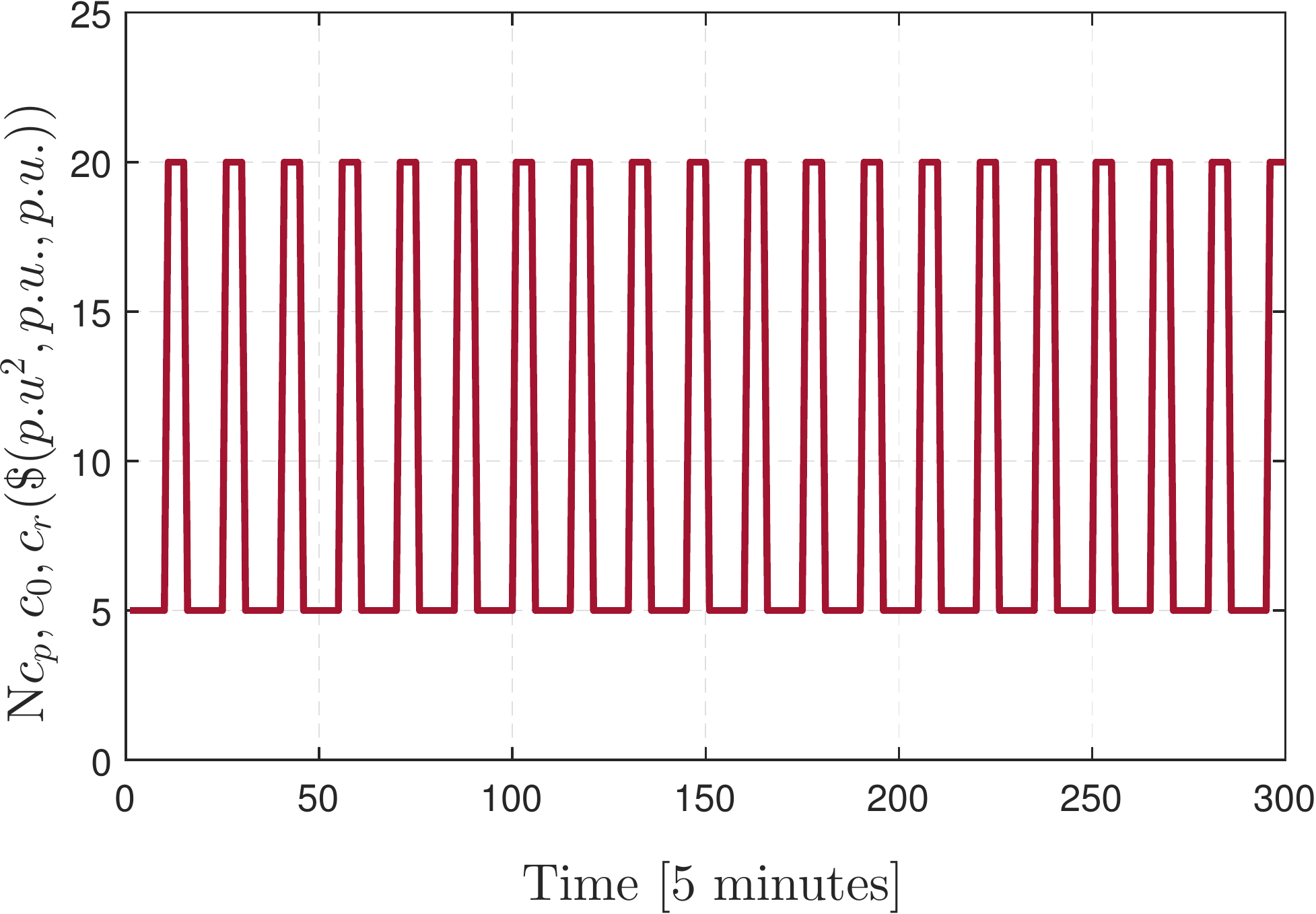}
\caption{Synthetic regulation signal $r^t$ (top) and prices $\{c^t_0,c^t_p,c^t_r\}$ (bottom) for Scenario 1.}\label{fig:case1}
\end{figure}

\begin{figure}[t]
\centering
\includegraphics[scale=0.3]{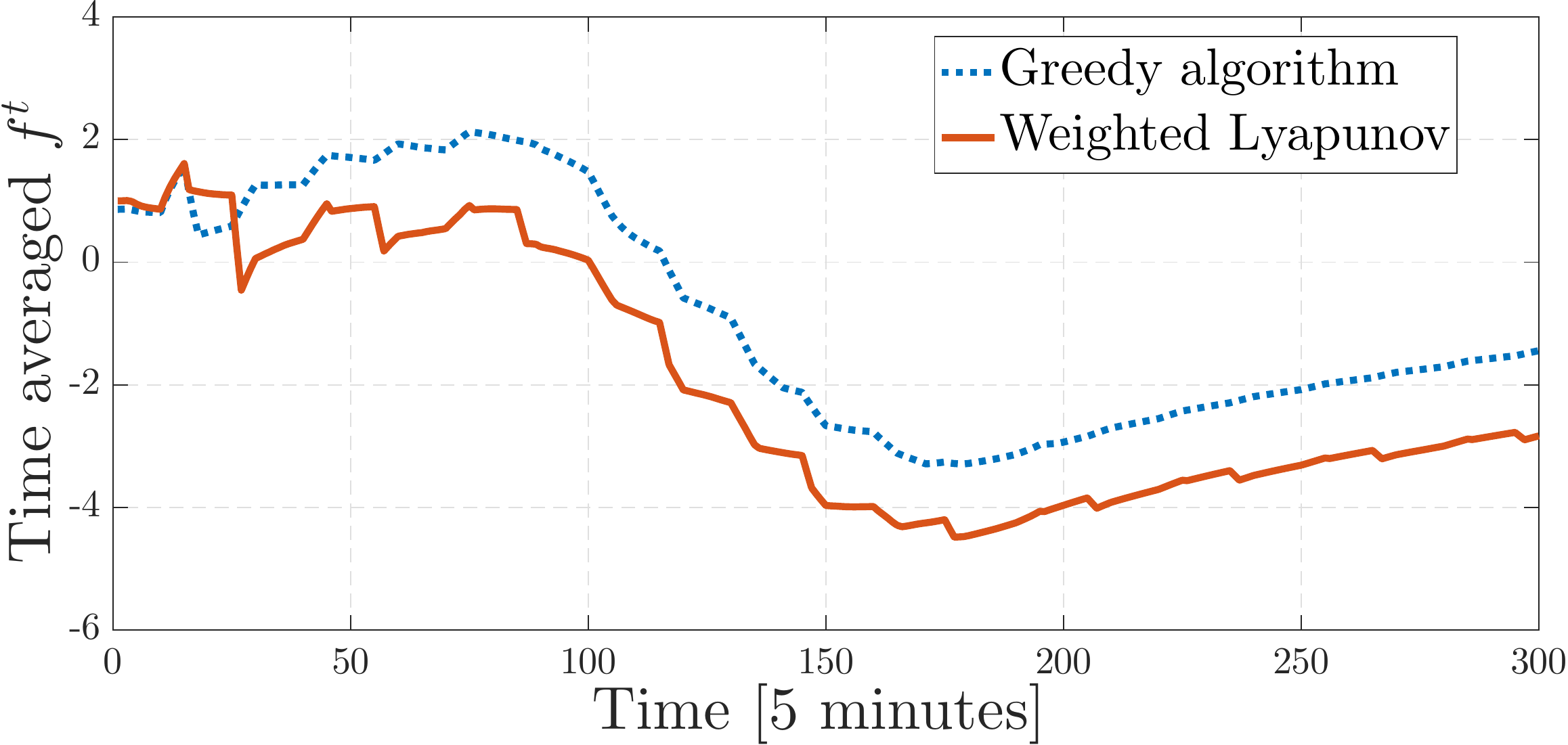}	%\vspace*{-5mm}
\caption{Comparison of average costs for Scenario 1.}
\label{fig:resultcase1}
\end{figure}

To showcase the superiority of the Lyapunov scheme over the greedy approach, we first tested the costs attained for the synthetically generated pricing and regulation signals shown in Fig.~\ref{fig:case1}. According to this setup named Scenario 1, the values for $r^t$ were oscillating between $\{\pm 1\}$ every 15 slots, and the prices $\{Nc^t_p, c^t_0,c_r^t\}$ were oscillating between $\{5,20\}~\$/\text{unit}$ with the former lasting for 10 slots and the latter for 5. Homes from the Pecan Street project with data identifiers 93, 171, 187, 252, 370, 545, 555, 585, 624, 744, 861, and 890 were placed on the buses of the IEEE 13-bus feeder of Fig.~\ref{fig:ieee13}. Further, we set $\underline{b}_n=-\overline{b}_n$ and $\underline{s}_n=0$ for all $n$. The average aggregate costs attained are depicted in Figure~\ref{fig:resultcase1}, where the weighted Lyapunov-based scheme clearly outperforms the greedy one. This is because the greedy scheme (dis)-charges the energy storage units myopically to their capacities during the low price of $\$5/\text{unit}$, rather than waiting to reap maximum rewards at $\$20/\text{unit}$. The Lyapunov scheme on the other hand saves some storage capacity for later opportunities.

To simulate a more realistic setup termed Scenario 2 was tested. Under Scenario 2, the price $c^t_0$ was set to the hourly real-time locational marginal prices for the RTO hub in the PJM market for 2011. Hourly prices were repeated 12 times to yield 5-minute prices. The coefficient $c^t_p$ was selected as $c^t_p=\frac{c^t_0}{N}$. Similarly, the price $c^t_r$ was set to the PJM regulation market clearing price for the same year, while the regulation signal $r^t$ was modeled as a zero-mean $\{\pm 1\}$ Bernoulli random variable capturing the nature of actual frequency regulation signals. The (dis)-charging rates were decreased by a factor of 16 compared to Scenario 1. 

\begin{figure}[t]
	\centering
	\includegraphics[scale=0.32]{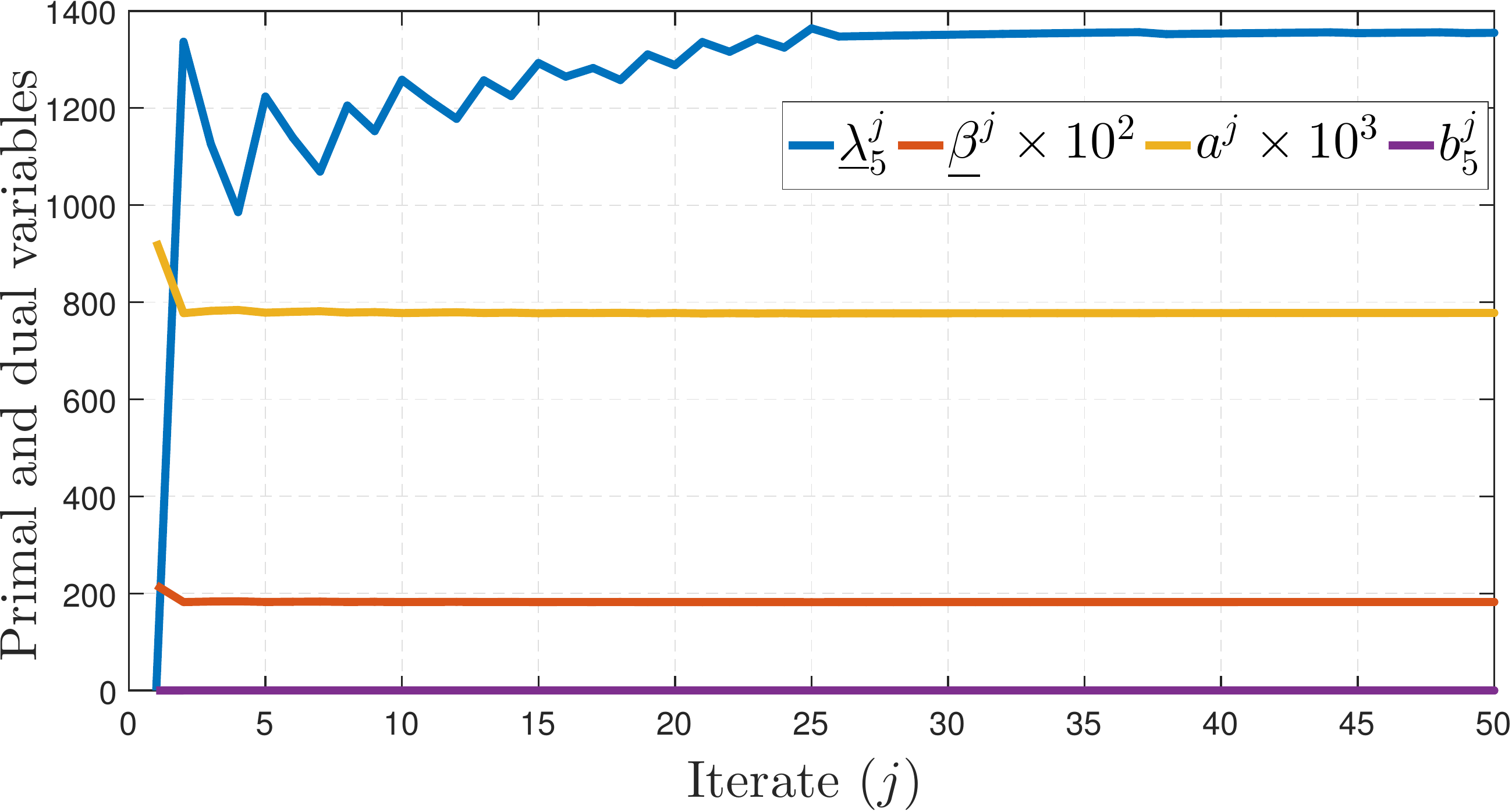}
	\caption{Convergence of primal and dual variables for bus 5.}
	\label{fig:primal}
\end{figure}

To demonstrate the convergence of Algorithm~\ref{alg:1}, Figure~\ref{fig:primal} shows the primal and dual variables corresponding to bus 5. During this period, the under-voltage constraint of \eqref{eq:voltreg} was active, thus yielding $\underline{\lambda}_5>0$. Using the diminishing step-size sequences $\eta_\nu^j=3\cdot 10^6/(j+1)$ and $\eta_\nu^j=2/(j+1)$, convergence was achieved within 30 iterations.

\begin{figure}[t]
	\centering
	\includegraphics[scale=0.32]{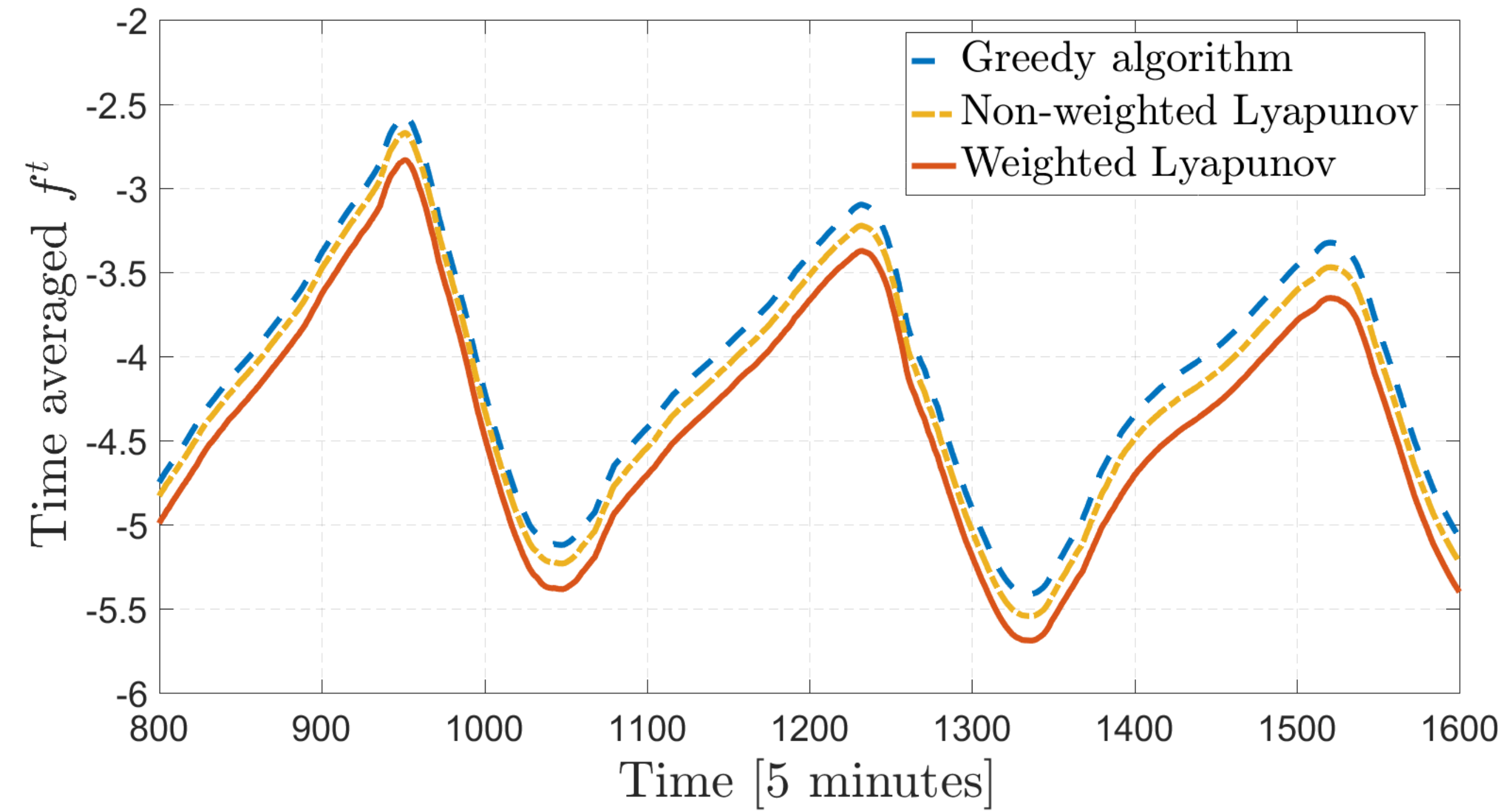}
	\caption{Averaged $f^t$ attained by different charging schemes for Scenario 2.}
	\label{fig:avg_cost}
\end{figure}

\begin{figure}[t]
	\centering
	\hspace*{-1em}
	\includegraphics[scale=0.34]{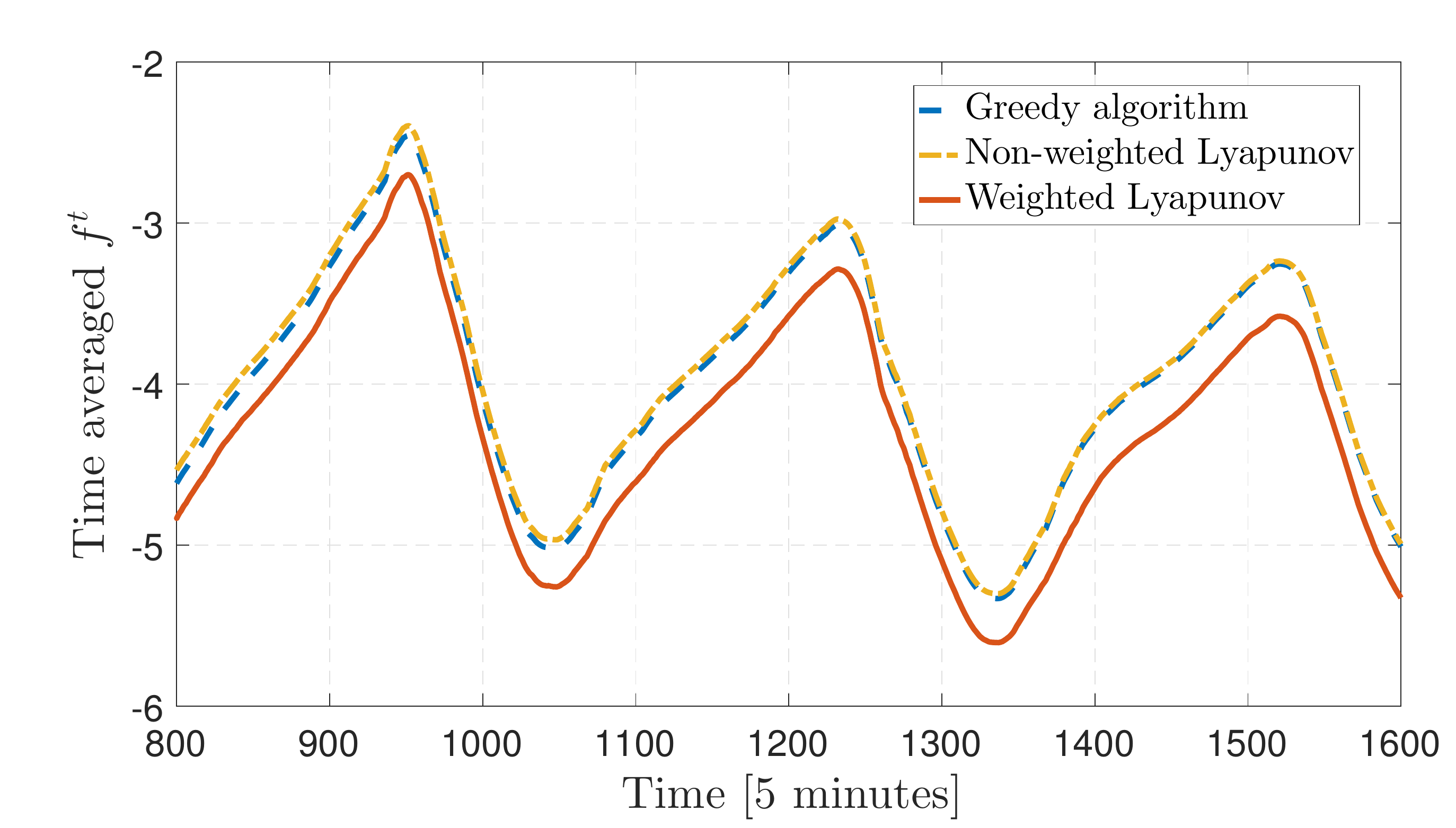}
	\caption{Averaged $f^t$ attained for Scenario 2 upon relaxing \eqref{eq:sign}.}
	\label{fig:avg_cost_noconst}
\end{figure}

Figure~\ref{fig:avg_cost} shows the time-averaged cost of \eqref{eq:centralized} obtained by the different schemes. The results indicate that the developed real-time solver outperforms both alternatives. Its superiority over \eqref{eq:greedy} is attributed to the myopic nature of the greedy scheme as discussed earlier. The improvement of our scheme over the standard Lyapunov approach is explained by the enhanced suboptimality gap of \eqref{eq:Kprime}. The larger suboptimality gap of the non-weighted scheme has been explained in the paragraph after Remark~\ref{co:maxopt}. Scenario 2 was also tested under the setup of Remark~\ref{re:no-r}, where charging decisions do not have to align with the regulation signal $r^t$. Even after this modification, approximately $85\%$ of the charging decisions still aligned with \eqref{eq:sign}, and the SOCs always respected \eqref{eq:soc}. The developed real-time solver still outperformed its alternatives as validated in Fig.~\ref{fig:avg_cost_noconst}.

\begin{figure}[t]
\centering
\includegraphics[scale=.3]{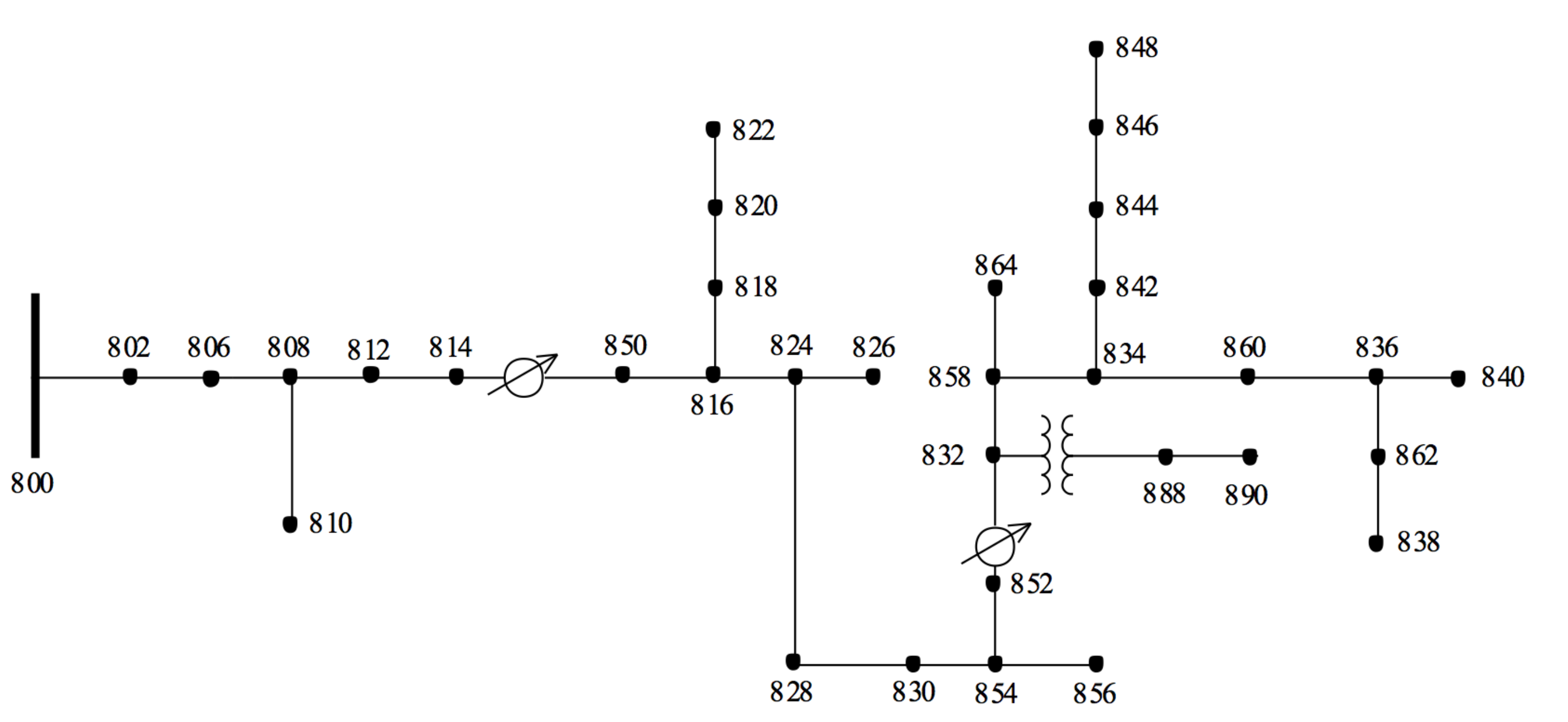}
\caption{IEEE 34-bus feeder for Scenario 3.}
%\vspace*{-1em}
\label{fig:ieee34}
\end{figure}

\begin{table}[t]
{{
\renewcommand{\arraystretch}{1.1}
\caption{Placing Load Data on the IEEE 34-bus Feeder.}\vspace*{-0.5em}
\label{tbl:Scenario3} \centering
\begin{tabular}{|rrr|rrr|}
\hline\hline
Bus \#\!&\! Home \# & $\overline{s}_n~[10^{-2}]$ & Bus \#\!&\! Home \# & $\overline{s}_n~[10^{-2}]$\\
\hline
\hline $802$ & $93$ &  $15$ & $806$ & $171$ & $7.7$ \\
\hline $808$ & $187$ & $7.9$ & $810$ & $252$ & $7.5$\\
\hline $812$ & $370$ & $11$ &  $814$ & $545$ & $0.5$\\
\hline $816$ & $744$ & $8.3$ & $818$ & $890$ & $6.1$\\
\hline $820$ & $1185$ & $15$ & $822$ & $1642$ & $7.9$\\
\hline $824$ & $861$ & $6.3$ & $826$& $1169$ & $6$\\
\hline $828$ & $1103$ & $11$ & $830$ & $1464$ & $7.7$\\
\hline $832$ & $3961$ & $8$ & $834$ & $6941$ & $11$\\
\hline $836$ & $8597$ & $7.9$ & $838$ & $9019$ & $8.3$\\
\hline $840$ & $8419$ & $5$ & $842$ & $8419$ & $7.7$\\
\hline $844$ & $9019$ & $7.5$ & $846$ & $8597$ & $8$\\
\hline $848$ & $9982$ & $6.3$ & $850$ & $624$ & $8$\\
\hline $852$ & $2980$ & $5$ & $854$ & $1718$ & $7.5$\\
\hline $856$ & $2129$ & $11$ & $858$ & $5129$ & $6.3$\\
\hline $860$ & $8084$ & $15$ & $862$ & $9982$ & $11$\\
\hline $864$& $6990$ & $6$ & $888$ & $4447$ & $8.3$\\
\hline $890$ & $5615$ & $6.1$ & & &\\
\hline \hline
\end{tabular}
}}
\end{table}

\begin{figure}[t]
\centering
\includegraphics[scale=0.33]{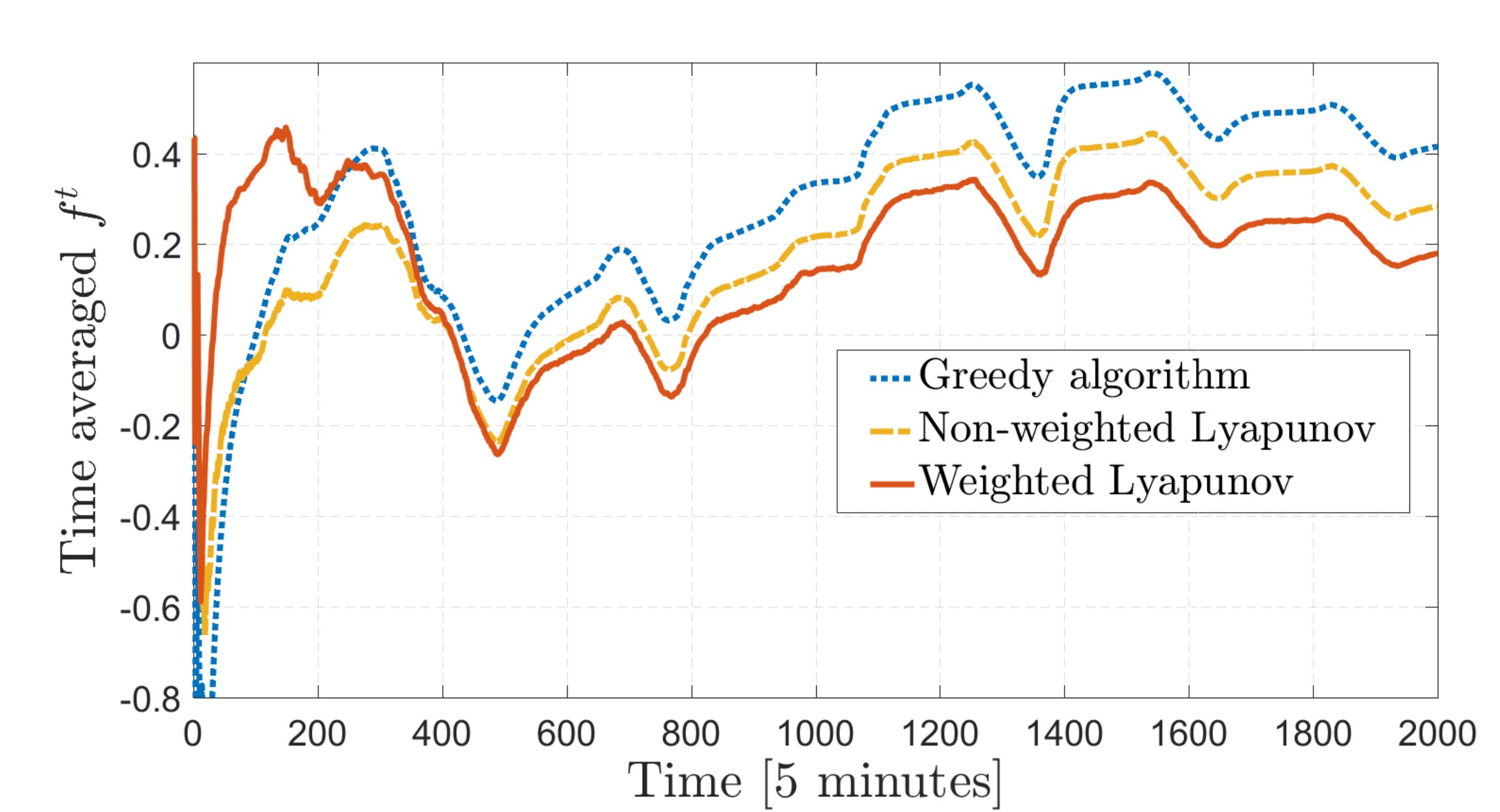}
\caption{Averaged $f^t$ attained by different charging schemes for Scenario 3.}
\label{fig:avg_cost2}
\end{figure} 

Finally, to study its scalability, the proposed scheme was tested on the IEEE 34-bus feeder shown in Fig.~\ref{fig:ieee34}. For this setup named Scenario 3, load data from the Pecan Street project were mapped to the feeder buses according to Table~\ref{tbl:Scenario3}. The energy storage parameters were set as $\overline{b}_n=\overline{s}_n/100$, $\underline{b}_n=-\overline{b}_n$, and $\underline{s}_n=0$ for all $n$. Tests were carried out using MATLAB R2016a on a 64-bit Windows 10 PC powered by a 2.6 GHz Intel i7-6700HQ CPU and 12 GB DDR3 RAM. The attained averaged costs are shown in Fig.~\ref{fig:avg_cost2}. The distributed implementation of Section~\ref{sec:dist} was timed against solving \eqref{eq:realtime} using YALMIP and the SeDuMi solver~\cite{YALMIP,sedumi}. The off-the-shelf solver ran for 430~sec, while the distributed algorithm needed only 10~sec.

%%%%%%%%%%%%%%%%%%%%%%%%%%%%%%%%%%%%%%%%%%%%%%%%%%%%%%%%%%%%%%%%%%	
\section{Conclusions}\label{sec:conclusions}
A novel approach combining game theory with Lyapunov optimization has been put forth to analyze competitive energy storage problems. A generalized Nash equilibrium has been shown to exist and can be found through a potential function. Leveraging Lyapunov optimization, a real-time scheme offering feasible charging decisions with suboptimality guarantees has been devised. The suggested decentralized implementation utilizes the distribution grid response and protects user's privacy. Numerical tests using realistic datasets have demonstrated the convergence of the distributed solver and the performance gain of the real-time scheme over its non-weighted counterpart and a greedy alternative. Extending our solvers to exact grid models and demand-response setups; considering multi-phase networks; including thermostatically-controlled loads; and incorporating partial information on future loads and prices form pertinent open research topics. 
 
\appendix
%%% Proof of Theorem characterizing the solution  %%%%%
\begin{IEEEproof}[Proof of Theorem~\ref{th:charac}] Arguing by contradiction, assume that hypothesis (a) holds for user $n$ yet $\hat{b}_n^t>0$. The case $\hat{b}_n^t<0$ is excluded because $r^t>0$ is assumed. Construct vector $\check{\mathbf{b}}^t$ with $\check{b}_n^t=0$ and $\check{b}^t_i=\hat{b}^t_i$ for all $i\neq n$. Under Assumption (a2) and because matrix $\mathbf{R}$ has non-negative entries, if $\hat{\mathbf{b}}^t$ is feasible for \eqref{eq:realtime}, then $\check{\mathbf{b}}^t$ is also feasible. 

It will be next shown that $\check{\mathbf{b}}^t$ attains an objective value for \eqref{eq:realtime} smaller than or equal to the one attained by $\hat{\mathbf{b}}^t$, that is
\begin{equation}\label{eq:costdiff}
\sum_{n=1}^N w_n x_n^t \check{b}_n^t + f^t(\check{\mathbf{b}}^t)\leq \sum_{n=1}^N w_n x_n^t \hat{b}_n^t + f^t(\hat{\mathbf{b}}^t).
\end{equation}
To this end, the difference of the upper bound values attained by $\hat{\mathbf{b}}^t$ and $\check{\mathbf{b}}^t$ is $\sum_{n=1}^N w_n x_n^t (\hat{b}_n^t - \check{b}_n^t) =\hat{b}_n^t w_n x_n^t $, while the difference of the instantaneous costs can be shown to be
\begin{equation*}%\label{eq:fdiff}
f^t(\hat{\mathbf{b}}^t)-f^t(\check{\mathbf{b}}^t)=\hat{b}_n^t\left[c^t_0+c^t_p\sum_{i=1}^N(\hat{b}^t_i+\ell^t_i)+c^t_p\ell^t_n- c_r^t\right].
\end{equation*}
Given $\hat{b}_n^t>0$, the inequality in \eqref{eq:costdiff} holds only if
\begin{equation}\label{eq:cond}
w_n x^t_n+c^t_0 + c^t_p\sum_{i=1}^{N}(\hat{b}^t_i+\ell^t_i)+c^t_p\ell^t_n- c^t_r \geq 0.
\end{equation}
Observe that $\hat{\mathbf{b}}^t\geq \mathbf{0}$ since $r^t>0$; loads are lower bounded by $\boldsymbol{\ell}^t\geq \underline{\boldsymbol{\ell}}$ from \eqref{eq:llim}; and prices $c_p^t$ and $c_r^t$ are bounded too. Then, the minimum value for $c^t_0+c^t_p\sum_{i=1}^{N}(\hat{b}^t_i+\ell^t_i)+c^t_p\ell^t_n- c^t_r$ is $\underline{g}_n$ by definition, and the hypothesis in (a) implies that \eqref{eq:cond} holds true. It has been shown that $\check{\mathbf{b}}^t$ yields the same or a lower cost than the unique minimizer $\hat{\mathbf{b}}^t$ of \eqref{eq:realtime} does. The latter is a contradiction and proves the claim. Claim (b) can be shown in a similar fashion.
\end{IEEEproof}

%%% Proof of Theorem on feasibility  %%%%%
\begin{IEEEproof}[Proof of Theorem~\ref{th:feasibility}] Proving by induction across time $t$, the base case holds true since $\hat{\mathbf{s}}^0\in [\underline{\mathbf{s}},\overline{\mathbf{s}}]$. Assuming $\hat{\mathbf{s}}^t\in [\underline{\mathbf{s}},\overline{\mathbf{s}}]$, it will be ensured that $\hat{\mathbf{s}}^{t+1}\in [\underline{\mathbf{s}},\overline{\mathbf{s}}]$. The analysis is performed on a per-user basis and over three cases:

\textbf{Case 1:} $x^t_n+\frac{\underline{g}_n}{w_n} \geq 0$. Depending on the regulation signal $r^t$, two subcases are considered. If $r^t>0$, then Th.~\ref{th:charac} asserts that $\hat{b}^t_n=0$ and $\hat{s}^{t+1}_n=\hat{s}^t_n$ so that the state remains feasible. 

If $r^t<0$, then $\hat{b}^t_n \in [\ub_n,0]$ and only the lower limit on $\hat{s}^{t+1}_n$ has to be ensured. To that end, substitute \eqref{eq:xs} in the assumption for Case 1 to get $\hat{s}^t_n +\gamma_n +\frac{\underline{g}_n}{w_n} \geq 0$. Combining the last inequality with the state update $\hat{s}^{t+1}_n = \hat{s}^t_n+\hat{b}^t_n$ yields:
\begin{equation}\label{eq:st1}
\hat{s}^{t+1}_n\geq -\gamma_n - \frac{\underline{g}_n}{w_n} +\hat{b}^t_n.
\end{equation}
The lower limit on $\hat{s}^{t+1}_n$ will be satisfied if the minimum value of the RHS in \eqref{eq:st1} is larger or equal to $\us_n$. Since $\hat{b}^t_n \in [\ub_n,0]$, the latter is guaranteed if 
\begin{equation}\label{eq:gammalim1}
\gamma_n\leq- \frac{\underline{g}_n}{w_n} +\ub_n- \us_n.
\end{equation}
	
\textbf{Case 2:} $x^t_n+\frac{\overline{g}_n}{w_n} \leq 0$. If $r^t<0$, then Th.~\ref{th:charac} guarantees $\hat{b}^t_n=0$ and the updated state $\hat{s}^{t+1}_n=\hat{s}_n^{t}$ remains feasible.

If $r^t>0$, then $\hat{b}^t_n \in [0,\ob_n]$ and only the upper limit on $\hat{s}^{t+1}_n$ needs to be ensured. From \eqref{eq:xs} and the assumption of Case 2, it follows that $\hat{s}^t_n+ \gamma_n +\frac{\overline{g}_n}{w_n}\leq 0$. Combining the last inequality with the state update $\hat{s}^{t+1}_n = \hat{s}^t_n+\hat{b}^t_n$ yields:
\begin{equation}\label{eq:st2}
\hat{s}^{t+1}_n \leq -\gamma_n - \frac{\overline{g}_n}{w_n} +\hat{b}^t_n.
\end{equation}
The upper limit on $\hat{s}^{t+1}_n$ will be satisfied if the maximum value of the RHS in \eqref{eq:st2} is smaller or equal to $\os_n$. Since $\hat{b}^t_n \in [0,\ob_n]$, the latter is guaranteed if 
\begin{equation}\label{eq:gammalim2}
\gamma_n\geq- \frac{\overline{g}_n}{w_n} +\ob_n- \os_n.
\end{equation}

\textbf{Case 3:} $-\frac{\overline{g}_n}{w_n}\leq x^t_n \leq -\frac{\underline{g}_n}{w_n}$. If $r^t>0$ then $\hat{b}^t_n \in [0,\ob_n]$ and only the upper limit on $\hat{s}^{t+1}_n$ needs to be maintained. Substituting \eqref{eq:xs} in the second inequality of Case 3 provides $\hat{s}^t_n +\gamma_n\leq -\frac{\underline{g}_n}{w_n}$. Substituting the state transition into the last inequality yields:
\begin{equation}\label{eq:st3}
\hat{s}^{t+1}_n \leq -\gamma_n - \frac{\underline{g}_n}{w_n} + \hat{b}^t_n
\end{equation}
The upper limit on $\hat{s}^{t+1}_n$ is respected if the maximum value of the RHS in \eqref{eq:st3} is $\le\os_n$. Since $\hat{b}^t_n \in [0,\ob_n]$, the latter is guaranteed if 
\begin{equation}\label{eq:gammalim3}
\gamma_n \geq - \frac{\underline{g}_n}{w_n} +\ob_n -\os_n.
\end{equation}
Because $\overline{g}_n > \underline{g}_n$, the bound of \eqref{eq:gammalim3} is tighter than the one in \eqref{eq:gammalim2}, and thus \eqref{eq:gammalim3} provides the lower bound on $\gamma_n$ in \eqref{eq:gamma}.

If $r^t<0$, then $\hat{b}^t_n \in [\ub_n,0]$ and only the lower limit on $\hat{s}^{t+1}_n$ needs to be maintained. Substituting \eqref{eq:xs} in the first inequality of Case 3 provides $-\frac{\overline{g}_n}{w_n} \leq \hat{s}^t_n + \gamma_n$. Substituting the state transition into the last inequality yields
\begin{equation}\label{eq:st4}
\hat{s}^{t+1}_n\geq -\gamma_n -\frac{\overline{g}_n}{w_n} +\hat{b}^t_n.
\end{equation}
Symmetrically to \eqref{eq:st3}, the lower limit on $\hat{s}^{t+1}_n$ is respected if the minimum value of the RHS in \eqref{eq:st4} is larger or equal to $\us_n$. Since $\hat{b}^t_n \in [\ub_n,0]$, the latter is guaranteed if 
\begin{equation}\label{eq:gammalim4}
\gamma_n\leq- \frac{\overline{g}_n}{w_n} +\ub_n- \us_n.
\end{equation}
Since the bound of \eqref{eq:gammalim4} is tighter than the one in \eqref{eq:gammalim1}, it is the former that determines the upper bound on $\gamma_n$ in \eqref{eq:gamma}.

Finally, it is easy see that the condition in \eqref{eq:mulim} guarantees that the bounds on $\gamma_n$ in \eqref{eq:gamma} yield a non-empty interval for all $n\in \mathcal{N}$.
\end{IEEEproof}

%%%%%%%%%%%%%%% Proof of Theorem on cost suboptimality bound  %%%%%%%%%%%%%%%%%%%
\begin{IEEEproof}[Proof of Theorem~\ref{th:subopt}]
Using \eqref{eq:driftpenalty}, the drift plus penalty term can be upper bounded as $\Delta^t+\mathbb{E}\left[f^t(\mathbf{b}^t)|\mathbf{x}^t\right] \leq\mathbb{E}\left[\sum_{n=1}^{N}w_n x^t_n b^t_n +  f^t(\mathbf{b}^t)\vert \mathbf{x}^t\right] +\frac{1}{2}\sum_{n=1}^{N}w_n\max\{\ob_n^2,\ub_n^2\}$. Note that the sequence of charging decisions $\{\hat{\mathbf{b}}^t\}$ obtained by \eqref{eq:realtime} essentially minimizes the aforementioned bound. Hence, the value attained for this bound by $\{\hat{\mathbf{b}}^t\}$ would be the minimum over all feasible policies, including the stationary policy $\{\grave{\mathbf{b}}^t\}$ of Lemma~\ref{th:stat}. Then, it follows that
\begin{equation}\label{eq:bound1}
\Delta^t+\mathbb{E}\left[f^t(\hat{\mathbf{b}}^t)|\mathbf{x}^t\right]\leq \phi^\star+\frac{1}{2}\sum_{n=1}^{N}w_n\max\{\ob_n^2,\ub_n^2\}.
\end{equation} 
Summing \eqref{eq:bound1} over $t=1,\ldots,T$; substituting $\Delta^t:=\mathbb{E} [L^{t+1}-L^t|\mathbf{x}^t ]$; and applying the law of total expectation yields
\begin{equation*}
\mathbb{E} [L^T - L^0] + \sum_{t=0}^{T-1}\mathbb{E} [f^t(\hat{\mathbf{b}}^t)]\leq  T\phi^\star + \frac{T}{2}\sum_{n=1}^{N} w_n \max\{\ob_n^2,\ub_n^2\}.
\end{equation*}
Because $\mathbb{E}[L^T]\geq 0$, the previous inequality yields
\begin{equation*}
\sum_{t=0}^{T-1}\mathbb{E} [f^t(\hat{\mathbf{b}}^t)]\leq  T\phi^\star + \frac{T}{2}\sum_{n=1}^{N} w_n \max\{\ob_n^2,\ub_n^2\} + \mathbb{E} [L^0].
\end{equation*}
Since $\mathbb{E}[L^0]$ is finite, dividing both sides by $ T$ and taking the limit of $T$ to infinity proves the claim.
\end{IEEEproof}

%%%%%%%%%%%%%%%%%% Proof of Theorem on cost minimizer bound  %%%%%%%%%%%%%%%%%%%
\begin{IEEEproof}[Proof of Theorem~\ref{th:l2norm}] The function $f^t$ is convex quadratic with Hessian matrix $\mathbf{H}^t:=c_p^t (\mathbf{I}+\mathbf{1}\mathbf{1}^\top)$. The minimum eigenvalue of $\mathbf{H}^t$ is $\lambda_{\min}(\mathbf{H}^t)=c_p^t$. Exploiting the Rayleigh quotient property of $\lambda_{\min}(\mathbf{H}^t)$ and using a second-order Taylor's series expansion of $f^t$ provides
\begin{equation}\label{eq:taylor}
f^t(\hat{\mathbf{b}}^t)\geq f^t(\tilde{\mathbf{b}}^t) + (\hat{\mathbf{b}}^t-\tilde{\mathbf{b}}^t)^\top \nabla f^t (\tilde{\mathbf{b}}^t)+\frac{c_p^t}{2}\|\hat{\mathbf{b}}^t-\tilde{\mathbf{b}}^t\|^2_2
\end{equation}
for all $t$. Applying the expectation operator over the random variables $\{c^t_{p},c^t_{r},\boldsymbol{\ell}^t,\mathbf{q}^t\}$; averaging over $t=1,\ldots,T$; and taking the limit of $T$ to infinity yields
\begin{align}
\lim_{T\to\infty} \frac{1}{T}\sum_{t=0}^{T-1}&\mathbb{E} [f^t(\hat{\mathbf{b}}^t)] \geq  \lim_{T\to\infty} \frac{1}{T} \sum_{t=0}^{T-1}\mathbb{E}[f^t({\tilde{\mathbf{b}}}^t)]\nonumber\\
&+\lim_{T\to\infty} \frac{1}{T}\sum_{t=0}^{T-1}\mathbb{E}[(\hat{\mathbf{b}}^t-\tilde{\mathbf{b}}^t)^\top \nabla f^t(\tilde{\mathbf{b}}^t)]\nonumber\\
&+\frac{\underline{c}_p}{2}\lim_{T\to\infty} \frac{1}{T} \sum_{t=0}^{T-1}\mathbb{E} [\|\hat{\mathbf{b}}^t-\tilde{\mathbf{b}}^t\|^2_2].\label{eq:taylor}
\end{align}
By the first-order optimality conditions for $\{\tilde{\mathbf{b}}^t\}$, it holds that
\begin{equation}
\lim_{T\to\infty} \frac{1}{T} \sum_{t=0}^{T-1}\mathbb{E} [ (\hat{\mathbf{b}}^t-\tilde{\mathbf{b}}^t)^\top \nabla f^t(\tilde{\mathbf{b}}^t) ]\geq 0. \label{eq:firstord}
\end{equation}
Note that $F(\{\mathbf{b}^t\})$ is a functional rather than a function. Nevertheless, it can be shown to be Fr\'{e}chet twice differentiable and optimality conditions for constrained optimization over functionals lead naturally to \eqref{eq:firstord}; see \cite[Prop.~2.11]{Burger}. Plugging \eqref{eq:firstord} and the definitions of $\tilde{\phi}$ and $\hat{\phi}$ into \eqref{eq:taylor} yields
\begin{equation*}
\lim_{T\to\infty} \frac{1}{T} \sum_{t=0}^{T-1}\mathbb{E} [\|\hat{\mathbf{b}}^t-\tilde{\mathbf{b}}^t\|^2_2 ]\leq \frac{2}{\underline{c}_p}(\hat{\phi}-\tilde{\phi}) \leq \frac{2K}{\underline{c}_p}
\end{equation*}
where the second inequality stems from \eqref{eq:subopt}.
\end{IEEEproof}

\bibliographystyle{IEEEtran}
\bibliography{myabrv,power}

\begin{IEEEbiography}[{\includegraphics[width=1in,height=1.25in,clip,keepaspectratio]{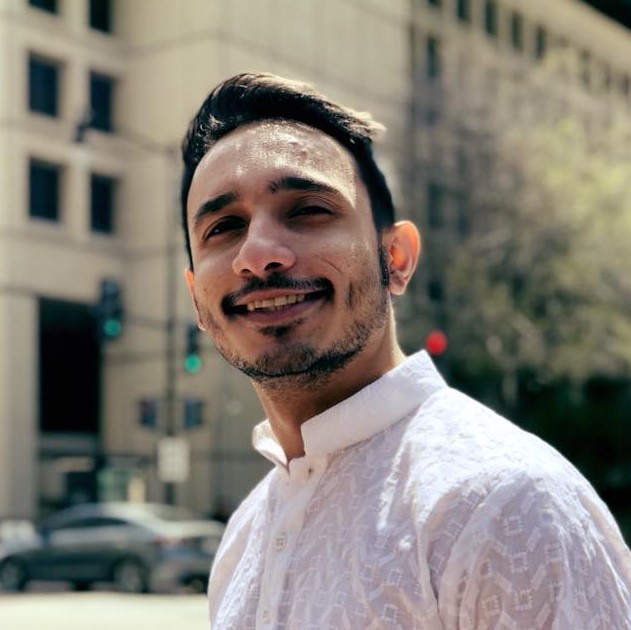}}]{Sarthak Gupta} received the B.Tech. degree in Electronics and Electrical Engnr. from IIT Guwahati, and the M.Sc. degree in power systems from Virginia Tech in 2013 and 2017, respectively. He is currently working as an Associate Engineer with the Distributed Resources Operations team at New York Independent System Operator (NYISO). His research interests include distributed resources, optimization algorithms, and renewable energy.
\end{IEEEbiography} 

\begin{IEEEbiography}[{\includegraphics[width=1in,height=1.25in,clip,keepaspectratio]{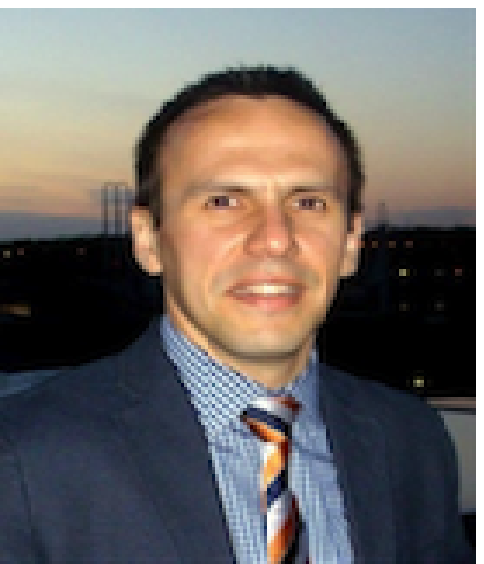}}] {Vassilis Kekatos} (SM'16) is an Assistant Professor of ECE at the Bradley Dept. of ECE at Virginia Tech. He obtained his Diploma, M.Sc., and Ph.D. in computer science and engr. from the Univ. of Patras, Greece, in 2001, 2003, and 2007, respectively. He is a recipient of the NSF Career Award in 2018 and the Marie Curie Fellowship during 2009-2012, and a research associate with the ECE Dept. at the Univ. of Minnesota, where he received the postdoctoral career development award (honorable mention). During 2014, he stayed with the Univ. of Texas at Austin and the Ohio State Univ. as a visiting researcher. His research focus is on optimization and learning for future energy systems. He is currently serving in the editorial board of the IEEE Trans. on Smart Grid.
\end{IEEEbiography}

\balance

\begin{IEEEbiography}[{\includegraphics[width=1in,height=1.25in,clip,keepaspectratio]{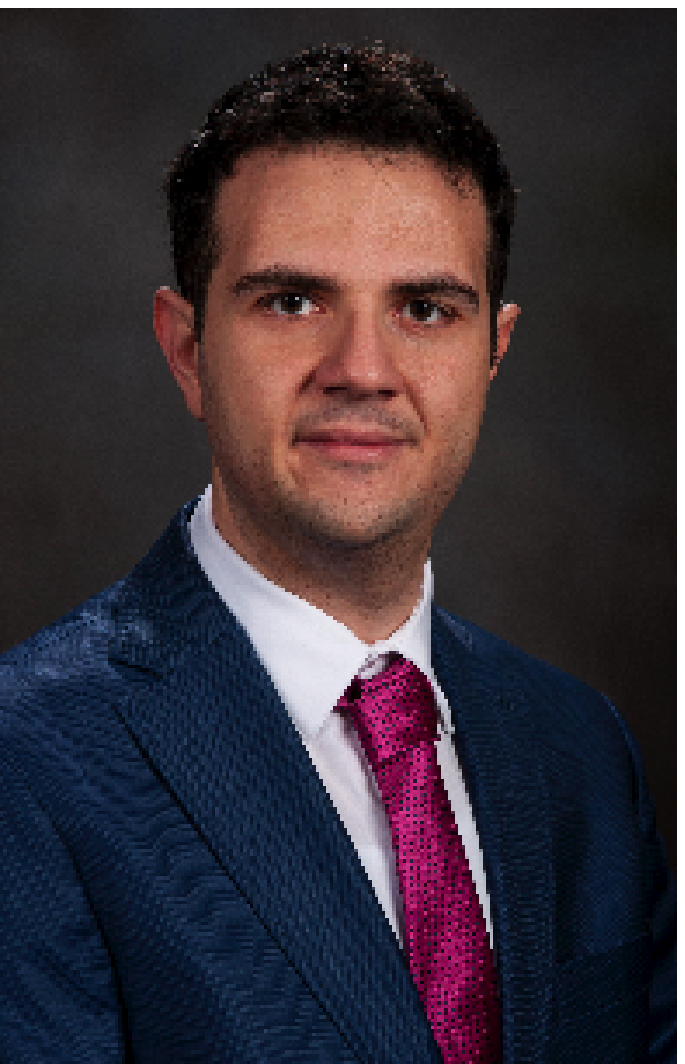}}] {Walid Saad} (S'07, M'10, SM’15) received his Ph.D degree from the University of Oslo in 2010. Currently,  he is an Associate Professor at the Department of Electrical and Computer Engineering at Virginia Tech, where he leads the Network Science, Wireless, and Security (NetSciWiS) laboratory, within the Wireless@VT research group. His  research interests include wireless networks, machine learning, game theory, cybersecurity, unmanned aerial vehicles, and cyber-physical systems. Dr. Saad is the recipient of the NSF CAREER award in 2013, the AFOSR summer faculty fellowship in 2014, and the Young Investigator Award from the Office of Naval Research (ONR) in 2015. He was the author/co-author of six conference best paper awards at WiOpt in 2009, ICIMP in 2010, IEEE WCNC in 2012,  IEEE PIMRC in 2015, IEEE SmartGridComm in 2015, and EuCNC in 2017. He is the recipient of the 2015 Fred W. Ellersick Prize from the IEEE Communications Society and of the 2017 IEEE ComSoc Best Young Professional in Academia award. From 2015-2017, Dr. Saad was named the Stephen O. Lane Junior Faculty Fellow at Virginia Tech and, in 2017, he was named College of Engineering Faculty Fellow. He currently serves as an editor for the IEEE Trans. on Wireless Communications, IEEE Trans. on Communications, IEEE Trans. on Mobile Computing, and IEEE Trans. on Information Forensics and Security.
\end{IEEEbiography}

\end{document}